\documentclass[journal,twoside,web]{ieeecolor}

\usepackage[utf8]{inputenc}
\usepackage{blindtext}
\usepackage{xcolor}
\usepackage{graphicx}
\usepackage{generic}
\usepackage{cite}
\usepackage{amsmath}

\usepackage{amsthm}
\usepackage{mathtools,amssymb,amsfonts}
\usepackage{graphicx}

\makeatletter
\let\NAT@parse\undefined
\makeatother
\usepackage{hyperref}  
\hypersetup{ colorlinks=true, linkcolor=blue, filecolor=magenta, urlcolor=cyan,}

\makeatletter
\renewcommand*\env@matrix[1][*\c@MaxMatrixCols c]{%
  \hskip -\arraycolsep
  \let\@ifnextchar\new@ifnextchar
  \array{#1}}
\makeatother

\usepackage{mathtools}
\usepackage{nccmath}
\usepackage{siunitx}
\usepackage{caption}
\usepackage{nccmath}
\usepackage{comment}

\usepackage{array}
\usepackage{multirow}

\usepackage{tikz}
\usetikzlibrary{fit,positioning,automata,backgrounds}

\newcommand{\matrixstyle}[1]{{#1}}

\newcommand{\expect}[1]{\mathbb{E}_{#1}}

\newcommand{\trace}{\mathrm{tr}}

\newtheorem{lemma}{Lemma}
\newtheorem{remark}{Remark}

\newtheorem{theorem}{Theorem}

\newtheorem{definition}{Definition}

\usepackage{enumerate}
\usepackage[ruled,vlined]{algorithm2e}
\usepackage{algorithmic}

\DeclareMathOperator*{\argmax}{arg\,max}  
\DeclareMathOperator*{\argmin}{arg\,min}  

\renewcommand{\baselinestretch}{0.975}

\def\BibTeX{{\rm B\kern-.05em{\sc i\kern-.025em b}\kern-.08em
    T\kern-.1667em\lower.7ex\hbox{E}\kern-.125emX}}
    
\begin{document}
\title{Deterministic Trajectory Optimization through Probabilistic Optimal Control}


\author{Mohammad Mahmoudi Filabadi, Tom Lefebvre and Guillaume Crevecoeur
\thanks{M. M. Filabadi, T. Lefebvre and G. Crevecoeur are with the Dynamic Design Lab (D2Lab) of the Department of Electromechanical, Systems and Metal Engineering, Ghent University, B-9052 Ghent, Belgium e-mail: 
{\tt\small \{mohammad.mahmoudifilabadi, tom.lefebvre, guillaume.crevecoeur\}@ugent.be}.}
\thanks{M. M. Filabadi, T. Lefebvre and G. Crevecoeur are members of the core lab MIRO, Flanders Make, Belgium.}
}

\maketitle


\begin{abstract}
In this article, we discuss two algorithms tailored to discrete-time deterministic finite-horizon nonlinear optimal control problems or so-called deterministic trajectory optimization problems. Both algorithms can be derived from an emerging theoretical paradigm that we refer to as probabilistic optimal control. The paradigm reformulates stochastic optimal control as an equivalent probabilistic inference problem and can be viewed as a generalisation of the former. The merit of this perspective is that it allows to address the problem using the Expectation-Maximization algorithm. It is shown that the application of this algorithm results in a fixed point iteration of probabilistic policies that converge to the deterministic optimal policy. Two strategies for policy evaluation are discussed, using state-of-the-art uncertainty quantification methods resulting into two distinct algorithms. The algorithms are structurally closest related to the differential dynamic programming algorithm and related methods that use sigma-point methods to avoid direct gradient evaluations. The main advantage of the algorithms is an improved balance between exploration and exploitation over the iterations, leading to improved numerical stability and accelerated convergence. These properties are demonstrated on different nonlinear systems.
\end{abstract}


\vspace{-10pt}

\section{Introduction} \label{sec:Intro}

Trajectory Optimization (TO) algorithms are a powerful class of methods for realizing goal-directed behaviour in dynamical systems. The aim is to find admissible state and action sequences by solving a discrete-time deterministic finite-horizon nonlinear optimal control problem. TO offers an elegant method both for offline and real-time control design in many application areas, such as robotics and autonomous vehicles. As a result (real-time) TO algorithms are a critical component of many contemporary control systems. Computation time is a critical factor driving algorithm selection. TO, especially for highly nonlinear and high-dimensional systems, thus remains an active research area in the control community, and lately also the reinforcement learning (RL) community.

A widely adopted method is Differential Dynamic Programming (DDP) \cite{Mayne1996,jacobson1970differential,tassa2014control}. The method updates an affine approximation of the policy backwards in time. This is done by approximation of the value function recursion, or Bellman equation (see section \ref{sec:differential-dynamic-programming}), with a quadratic surrogate. The solution is evaluated by forward simulation of the affine policy and the procedure is then iterated. Exact DDP is computationally relatively expensive as it requires computing the second-order derivatives of the dynamics. A popular simplification is the iterative Linear Quadratic Regulator (iLQR) \cite{Weiwei2004}, which uses only the first-order derivatives of the dynamics to reduce computation time, however, this comes at the cost of slower convergence. In recent work, the unscented transform and sparse Gauss-Hermite quadrature rules have been employed to build the quadratic surrogate, avoiding direct evaluation of any gradients \cite{Manchester2016,SGHQ_DDP_Shaoming2019}. The former are well-established numerical techniques used for uncertainty quantification and nonlinear filtering algorithms \cite{sarkka2023}. The Fourier-Hermite series and sigma points for approximating the action-value function were used in \cite{Hassan2023}. In the works mentioned, a Gaussian distribution is used to generate the sigma points. The covariance of the Gaussian is chosen arbitrarily and is constant over the time horizon. The authors emphasize the performance of their algorithm strongly depends on this choice. Automatically choosing the covariance remained an open problem \cite{Manchester2016,Hassan2023}. 

In this paper, we evaluate the probabilistic optimal control (POC) paradigm \cite{Karny1996,toussaint2006,toussaint2009robot,levine2018,rawlik2013stochastic,rawlik2013probabilistic,Noorani2022,lefebvre2023, watson2020stochastic, watson2021advancing}
 as a source of inspiration to derive or reimagine deterministic TO algorithms. As we will show, this will address the issue mentioned above, as well as offer other benefits. A key property of POC is that it re-establishes Stochastic Optimal Control (SOC) and Risk-Sensitive Optimal Control (RSOC) as probabilistic inference problems. If successful, indeed this means that the underlying optimal control problems could be treated by means of numerical inference methods. Particularly, here we are interested in the case where the stochastic dynamics reduce to deterministic dynamics, rendering the SOC problem into a deterministic TO problem. Then, the inference based methods could be used to treat the deterministic TO problem rather than the classical methods discussed above. Of course this presumes that the POC problem and its solution are equivalent to the underlying deterministic optimization problem and not merely similar. Ultimately, we do wish to solve the deterministic TO problem at hand.

Formal equivalence can be established  between the Risk-Sensitive Optimal Control (RSOC) problem and a particular Maximum Likelihood Estimation (MLE) problem as discussed in e.g. \cite{rawlik2013stochastic,rawlik2013probabilistic,Noorani2022,lefebvre2023}. The equivalence allows us to treat the RSOC problem using well-established inference techniques, in particular, the Expectation-Maximization (EM) algorithm \cite{bishop2006pattern}. The EM algorithm is a well-known strategy to treat MLE problems by replacing the direct and intractable optimization problem with a sequence of surrogate problems that are easier to treat. Application of the EM algorithm to the RSOC problem results into a fixed point iteration of probabilistic policies that converge to the deterministic RSOC policy in the limit \cite{lefebvre2023}. The uncertainty on the policy sequences can be interpreted as Bayesian or epistemic uncertainty on the deterministic optimal policy sequence. The uncertainty ultimately vanishes when the policy converges to the deterministic optimal policy.

This and similar ideas have a long history in the research community. We give here only a brief overview. A detailed and technical comparison with closely related work is given in a later section (section \ref{sec:prev_work}). The use of the EM algorithm to treat RL problems was discussed as early as \cite{Dayan1997}. The idea to design an approximate optimal control by matching a closed-loop density with some desired density was proposed about the same time \cite{Karny1996}. Here the desired density did not yet internalize the notion of an external cost or reward. A formal similarity between a particular probabilistic model and the SOC problem was first discussed in \cite{toussaint2006,toussaint2009robot}. Here the notion of an external cost was internalized into the probabilistic model used to represent the system therewith producing an optimal density. It was further proposed that the optimal control could be inferred through likelihood maximization. In the same work this was already brought into connection with Gaussian message-passing techniques or Bayesian Smoothing. Though it was also noted that likelihood maximization was not equivalent to minimization of the underlying SOC problem. In \cite{rawlik2013stochastic,rawlik2013probabilistic} it was remarked that maximizing the likelihood is in fact equivalent to solving an RSOC problem rather than a SOC problem. Application of the EM algorithm to address the maximum likelihood estimation (MLE) was not considered here explicitly. Rather the authors show that the associated optimal density, with internalized notion of cost, can be integrated in a density matching framework. This produced an entropy regularized SOC problem. Entropy regularization of SOC problems is a well-studied and -understood problem that was first discussed in \cite{ziebart2008maximum,ziebart2010modelingb}. There it was shown that the optimal policy was no longer deterministic but stochastic or probabilistic, i.e. is the optimal policy is a density. This idea was further explored to produce several efficient RL algorithms characterised by their ability to explore intelligently on account of the balanced entropy maximization of the policy
\cite{levine2013variational,levine2013guided,haarnoja2018soft}. In the same work that highlighted the equivalence between RSOC and MLE \cite{rawlik2013stochastic,rawlik2013probabilistic}, it was remarked that the optimal density contains a representation of a prior policy and that therefore an iterative procedure could be established where the optimal policy from one iteration is substituted as a prior in the next iteration. It was assumed that the resulting sequence would converge to optimal policy corresponding the associated RSOC problem. The authors further showed that the optimal policy corresponding the entropy regularized SOC could be evaluated by dynamic programming, and suggested, based on the analogy with the seminal work of \cite{toussaint2006,toussaint2009robot}, also by means of Bayesian Smoothing. 
Thereafter, researchers in \cite{watson2020stochastic,watson2021advancing,watson2021b} 
formulated the policy inference problem as an input estimation problem and framed it using an EM algorithm.
Notably, these studies offer the Bayesian Smoothing approach for stochastic control in the E-step and optimize a hyperparameter in the M-step of the EM algorithm. 
Additionally, the use of approximate Bayesian Smoothing approaches that leverage state-of-the-art uncertainty quantification methods, such as Gauss-Hermite and cubature quadrature, is discussed in \cite{watson2021b}.

Still there remains some ambiguity about what problem is solved exactly in the end, how the entropy regularized SOC problem relates exactly to the optimal density and the likelihood maximization problem from \cite{toussaint2006,toussaint2009robot}, etc. Some of these issues were discussed in \cite{levine2018}. Here it was also shown how to extract an optimal policy from the optimal density from \cite{toussaint2006,toussaint2009robot} by means of dynamic programming. Interestingly, the associated backward recursion is not the same as the backward recursion associated to the entropy regularized SOC problem and thus the associated optimal policies are also not equivalent. The author further discusses whether this optimal policy optimizes an interpretable objective however only reaches a satisfactory answer for the case of deterministic dynamics. 
In \cite{lefebvre2023} it was attempted to clarify several of these points. It was pointed out that iterating the solution to the entropy regularized SOC problem converges to the SOC policy, whereas solving the MLE problem associated to the RSOC problem converges to the RSOC policy. Only the latter can be treated using both dynamic programming as well as Bayesian smoothing or message passing. The former can be treated only using dynamic programming. It is noted that this difference disappears in the context of deterministic dynamics.

Based on these prior results, in this work, we pursue the following. First, we wish to provide a clear and encompassing development of the MLE treatment of the RSOC problem and how the solution can either be obtained through iterated dynamic programming or Bayesian smoothing. Our particular interest is to exploit the framework to address deterministic TO problems. To that end, insights are provided that support this interest. Secondly, we discuss how two derivative-free algorithms can be derived, which we refer to as the Sigma-Point Probabilistic Dynamic Programming (SP-PDP) and the Sigma-Point Bayesian Smoothing Control (SP-BSC) algorithm. Each algorithm pursues a different strategy to evaluate the iterate probabilistic policies, resulting in two distinct computational structures. 
Given the probabilistic nature of the approaches, we also rely on state-of-the-art uncertainty quantification methods to render the strategies into practical algorithms. The use of sigma points is directly or indirectly related to the evaluation of pseudo-gradients which can be loosely interpreted as probing the optimization landscape. The main difference when compared to other trajectory optimization algorithms such as Differential Dynamic Programming is that these algorithms maintain an approximation of the probabilistic policy rather than of the deterministic policy over the iterations and automatically tune a hyperparameter related to exploration. 

We organize the remainder of the paper as follows. Section \ref{sec:Preliminaries} states the problem formulation and preliminary background. In section \ref{sec:inference} we develop and treat the probabilistic reformulation of the nonlinear discrete-time deterministic optimal control problem. In section \ref{sec:algorithm} subsequently we discuss two algorithms by making use of the state-of-the-art numerical techniques in uncertainty quantification and propagation 
and discuss in detail the connection with previous methods. Numerical experiments are documented in section \ref{sec:results} on various nonlinear systems to show the capability of the algorithms in comparison to previous methods.  
The results demonstrate that the algorithms are able to handle nonlinear and high-dimensional systems similar to its predecessors.

\section{Background} \label{sec:Preliminaries} 
\subsection{Notation}
The set of non-negative real numbers is denoted $\mathbb{R}^+$ and the set of strictly positive real numbers is denoted $\mathbb{R}^+_*$. We write $\alpha_t$ to denote the value or instance of a variable or function, $\alpha$, at the discrete time instant $t$. Similarly, $\alpha_{t_1:t_2}$ refers to the collection of values or instances of a variable or function over the discrete time interval $[t_1,t_2]$, i.e. $\alpha_{t_1:t_2}=\{\alpha_{t_1},\dots,\alpha_{t_2}\}$. Further, we employ the asterisk, $*$, to denote symmetric entries in matrices. With $p(\cdot)$ we refer to any member of the class of probability density functions. The arguments of the corresponding function class are implied by the context.  Notation $p(A|B)$ denotes the probability density of $A$ conditional on $B$ whereas $p(A;B)$ denotes the probability density of $A$ parametrized by $B$. The expression $p(A|B;C)$ thus codifies the probability density of $A$ conditional on $B$ parametrized by $C$. Here $C$ can be a parameter or it can be a function. Finally, $\mathbb{E}[X]$ denotes the expected value of a random variable, $X$. When we subscript the operator this is to highlight the probability measure. 

\subsection{Problem formulation} \label{sec:problem_formulation}
We consider a class of nonlinear discrete-time deterministic systems governed by the following difference equation.
\begin{equation} \label{eq:general_dynamics_transition}
    x_{t+1} = f_t(x_t, u_t), ~t\in\{0,\dots,T-1\}
\end{equation}
Here $x_{t}\in \mathbb{R}^{n_x}$ and $u_{t}\in \mathbb{R}^{n_u}$ represent the state and action vector of the dynamic system at the discrete time $t$, and $f_t:\mathbb{R}^{n_x}\times\mathbb{R}^{n_u}\times\mathbb{R}\mapsto \mathbb{R}^{n_x}$ represents the state-space function. Henceforth, to denote the concatenated state-action vector at time $t$, we define $\xi_t \triangleq (x_t, u_t) \in \mathbb{R}^{n_{\xi}}$ with $n_{\xi} \triangleq n_x+n_u$.

Further, we consider the cost
\begin{equation} \label{eq:general_cost_functional}
    C(\xi_{0:T}) = c_T(x_T) + \sum\nolimits_{t=0}^{T-1} c_t(\xi_t)
\end{equation}
subject to the dynamics (\ref{eq:general_dynamics_transition}) with given initial state $x_0$. The stage cost at time $t$, is denoted, $c_t:\mathbb{R}^{n_\xi}\mapsto \mathbb{R}^+ \ge 0$, and, the terminal cost associated with the final state, $x_T$, is denoted $c_T:\mathbb{R}^{n_x}\mapsto \mathbb{R}^+$. For notational convenience we write $\xi_{0:T}$ when we mean $\{\xi_{0:T-1},x_T\}$. We will refer to this collection of variables as the state-action trajectory.

In this work, we aim to find the (deterministic) optimal control sequence, $u^*_{0:T-1}$, that minimizes the following nonlinear discrete-time deterministic optimal control problem.
\begin{equation} \label{eq:general_min_cost_functional}
    u^*_{0:T-1} = \argmin_{u_{0:T-1}} C(\xi_{0:T}).
\end{equation}
\subsection{Dynamic programming}\label{sec:dynamic-programming}
The solution to problem (\ref{eq:general_min_cost_functional}) can be represented in terms of the so-called state-action value function, $Q_t:\mathbb{R}^{n_\xi}\mapsto \mathbb{R}^+$ \cite{Bellman1957}. The policy map, $u_t^*:\mathbb{R}^{n_x}\mapsto \mathbb{R}^{n_u}$, determines the optimal action at time $t$ for given state $x_t$ and is referred to as the optimal policy. This representation of the solution is clearly more general than the representation in (\ref{eq:general_min_cost_functional}). The latter can be reproduced by a forward simulation of the nonlinear system dynamics (\ref{eq:general_dynamics_transition}) with the optimal policy in a closed loop.
\begin{equation}
	u_t^*(x_t) = \argmin_{u_t} Q_t(\xi_t)  \label{eq:optimal_determinstic}
\end{equation}

The state-action value function, $Q_t$, is defined as the superposition of the stage cost and the value function, $V_t:\mathbb{R}^{n_x}\mapsto\mathbb{R}^+$, or the so-called optimal cost-to-go. The value function is defined as the minimum cost that is accumulated between time step $t$ and $T$ starting from state $x_t$, when administering the optimal control in (\ref{eq:optimal_determinstic}) at every time step and present state. 
\begin{equation}
	Q_t(\xi_t) \triangleq c_t(\xi_t) + V_{t+1}(f_t(\xi_t))  \label{eq:Q_function} 
\end{equation}

The value function is governed by the backward recursive Bellman equation \cite{Bellman1957}, starting from the final time step $T$ to the present time step $t$. It is implied that $V_T(x_T) = c_T(x_T)$.
\begin{equation} \label{eq:Bellman_Optimality}
    V_t(x_t) = \min_{u_t} c_t(\xi_t) + V_{t+1}(f_t(\xi_t))
\end{equation}

Solving the backward recursion explicitly is computationally infeasible for general nonlinear problems. To address this issue in practice, one usually appeals to approximations of the state-action value function. 
One of the most popular and widely used algorithms for this purpose is DDP. The DDP algorithm relies on local quadratic approximation of the state-action value function. For completeness and future comparison, the DDP algorithm is covered in more detail next.

\subsection{Differential dynamic programming}\label{sec:differential-dynamic-programming}
The standard DDP algorithm \cite{Mayne1996, jacobson1970differential} assumes that the state-action value and optimal cost-to-go functions, $Q_{0:T-1}$ and $V_{0:T-1}$, can be accurately approximated by a quadratic surrogate about some nominal state-action trajectory, $\hat{\xi}_{0:T}$.
\begin{subequations} \label{eq:Func_Approx1}
	\begin{align} 
		 \label{eq:Func_Approx1_Q}
		Q_t({\xi}_t) \approx \hat{Q}_t(\xi_t) &= \tfrac{1}{2}\begin{bmatrix}
			1 \\ \delta {\xi}_{t} 
		\end{bmatrix}^\top \begin{bmatrix}
			2\hat{Q}_{0,t} & \hat{Q}_{\xi,t}^{\top}  \\
			* & \hat{Q}_{\xi\xi,t} \\
		\end{bmatrix}\begin{bmatrix}
			1 \\ \delta {\xi}_{t}
		\end{bmatrix} \\
			V_t({x}_t) \approx \hat{V}_t({x}_t) &= \tfrac{1}{2}\begin{bmatrix}
			1 \\ \delta {x}_{t} 
		\end{bmatrix}^\top \begin{bmatrix}
			2\hat{V}_{0,t} & \hat{V}_{x,t}^{\top} \\
			* & \hat{V}_{xx,t}
		\end{bmatrix}\begin{bmatrix}
			1 \\ \delta {x}_{t} 
		\end{bmatrix}
		 \label{eq:Func_Approx1_V}
	\end{align}
\end{subequations}
Here $\delta \xi_t \triangleq \xi_t - \hat{\xi_t}$ denote deviations from the nominal trajectory. The interpretation of $\delta {x}_t$ and $\delta {u}_t$ is trivial. Further, the following matrix partitions are implied.
\begin{subequations}
	\begin{align}
		\hat{Q}_{\xi,t} &\triangleq 
	\begin{bmatrix}
		\hat{Q}_{x,t}   \\
		\hat{Q}_{u,t}   \\
	\end{bmatrix} \\
	\hat{Q}_{\xi\xi,t} &\triangleq 
	\begin{bmatrix}
		\hat{Q}_{xx,t} & \hat{Q}_{ux,t}^{\top}  \\
		* & \hat{Q}_{uu,t} \\
	\end{bmatrix}
	\end{align}
\end{subequations}

The principle idea behind the DDP algorithm is then to evaluate the recursion in (\ref{eq:Bellman_Optimality}) by approximating (\ref{eq:Q_function}) using a second-order Taylor series expansion whilst also practising the assumption that $V_t$ is quadratic. One easily verifies the following expressions for the coefficients of $\hat{Q}_t$.
\begin{subequations} \label{eq:Q_update_DDP}
	\begin{align} 
		\hat{Q}_{0,t} &= \hat{c}_{0,t} + \hat{V}_{0,t+1}\\
		\hat{Q}_{\xi,t} &= \hat{C}_{\xi,t} + \hat{\matrixstyle{F}}_{\xi,t}^\top\hat{V}_{x,t+1} \\
		\hat{Q}_{\xi\xi,t} &= \hat{C}_{\xi\xi,t} + \hat{\matrixstyle{F}}_{\xi,t}^\top \hat{V}_{xx,t+1}\hat{\matrixstyle{F}}_{\xi,t} + \sum\nolimits_{i=1}^{n_x} \hat{\matrixstyle{F}}_{\xi \xi,t}^i \hat{V}_{x,t+1}^{i}
	\end{align}
\end{subequations}
Here the scalar $\hat{c}_{0,t}$, vector $\hat{C}_{\xi,t}$ and matrix $\hat{C}_{\xi\xi,t}$ represent the function value, Jacobian and Hessian of the stage cost, $c_t$, respectively, evaluated at $\hat{\xi}_t$. The matrices $\hat{F}_{\xi,t}$ and $\hat{F}_{\xi\xi,t}^i$ with $i\in\{1,\dots,n_x\}$ represent the Jacobian and Hessian matrices of the dynamic function, $f_t$, and, dynamic function elements, $f_t^i$, respectively, again, evaluated at $\hat{\xi}_t$.

Then, following (\ref{eq:optimal_determinstic}), we can obtain an optimal deviation, $\delta u_t^*$, on the nominal action, $\hat{u}_t$, by optimizing (\ref{eq:Func_Approx1_Q}) with respect to $\delta u_t$. A minimum is attained for
\begin{equation} \label{eq:optimal_cont_DDP}
    \begin{aligned}
        \delta u_t^* &= - \hat{Q}_{uu,t}^{-1} \hat{Q}_{u,t} - \hat{Q}_{uu,t}^{-1} \hat{Q}_{ux,t}\delta x_t, \\
    \end{aligned}
\end{equation}
Let us now also define the optimal control gains
\begin{subequations} \label{eq:DDP_opt_gain}
\begin{align}
	k^*_t &\triangleq -\hat{Q}_{uu,t}^{-1} \hat{Q}_{u,t}
	 \\ \matrixstyle{K}^*_t &\triangleq -\hat{Q}_{uu,t}^{-1} \hat{Q}_{ux,t} 
\end{align}
\end{subequations}

Finally, we can find explicit recursive expressions for the coefficients of $\hat{V}_t$ by substituting (\ref{eq:Func_Approx1_Q}), (\ref{eq:Func_Approx1_V}) and (\ref{eq:optimal_cont_DDP}) into (\ref{eq:Bellman_Optimality}). 
\begin{equation} \label{eq:value_updat_DDP}
				\begin{aligned} 
                    \hat{V}_{0,t} &= \hat{Q}_{0,t} - \tfrac{1}{2}k_t^{*\top} \hat{Q}_{uu,t} k^*_t \\
					\hat{V}_{x,t} &= \hat{Q}_{x,t} - \matrixstyle{K}_t^{*\top} \hat{Q}_{uu,t} k^*_t \\
					\hat{V}_{xx,t} &= \hat{Q}_{xx,t}  - \matrixstyle{K}_t^{*\top} \hat{Q}_{uu,t} \matrixstyle{K}^*_t 
				\end{aligned}
\end{equation}

The calculations above are at the heart of the DDP algorithm. The entire procedure then exists in practising a \textit{backward pass} where one constructs quadratic state-action value and optimal cost-to-go surrogates about the previous state-action trajectory, starting from a second-order Taylor expansion of $V_T(x_T) = c_T(x_T)$. Further one evaluates the recursive expressions in (\ref{eq:Q_update_DDP})-(\ref{eq:value_updat_DDP}) from time step $t = T-1$ until $t = 0$. In a \textit{forward pass} the system is then simulated forward in time administering the updated control input $u_t^* = \hat{u}_t + \delta u_t^*$ with $\delta x_t^* = x_t^* - \hat{x}_t$ which produces a new state-action trajectory, $\xi_{0:T}^*$, that will serve as the reference for another iteration and so forth. The backward and forward passes are thence iterated until convergence to obtain the desired optimal control trajectory, $u^*_{0:T-1}$, and associated state trajectory, $x^*_{0:T}$.

\begin{remark}
	The computational structure of the DDP algorithm is backward-then-forward. The policy is first calculated backwards in time and then evaluated forwards in time.
\end{remark}

\section{Probabilistic Optimal Control} \label{sec:inference}

The goal of this section is to cast and treat the problem (\ref{eq:general_min_cost_functional}) as an equivalent probabilistic inference problem. In the next section, we will discuss how algorithms can be derived from the general theory. To that end, the following steps are needed.

In sec. \ref{sec:problem-reformulation}, we first expand the solution space to the set of probabilistic policies. Formally, this does not affect the optimal solution yet this is an important step for the ensuing development. Second, we reformulate the optimal control problem as a Risk-Sensitive Optimal Control (RSOC) problem. This step is crucial to establish the equivalence. In sec. \ref{sec:equivalent-probabilistic-inference-problem}, we internalize the notion of reward/cost by introducing a set of dummy random variables. This will finally allow us to recast the RSOC problem as an equivalent Maximum Likelihood Estimation (MLE) problem. The generic treatment of the corresponding MLE problem using the Expectation-Maximization (EM) algorithm is then discussed in sec. \ref{sec:expectation-maximization}. Specific to the E-step is that it involves the evaluation of some posterior density from which optimal policies are derived in the M-step. In sec. \ref{sec:iterative-solution-methods} we discuss two procedures to evaluate this posterior and with that the optimal policies. {Interestingly, the first strategy is characterised by a \textit{backward-then-forward} computational structure whereas that of the second is \textit{forward-then-backward}. This manifests into different numerical algorithms in sec. \ref{sec:algorithm}.}





\subsection{Problem reformulation}\label{sec:problem-reformulation}

Although the dynamics of the original problem are deterministic, we can formally represent them probabilistically. It is straightforward to represent the dynamics of the system as a (controlled) Markov model. To that end, we define the following transition probability distribution, $p(x_{t+1}|\xi_t)$, and, the initial state distribution, $p(x_0)$. Here $\delta(\cdot)$ represents the Dirac delta distribution and $\mathfrak{x}_0$ denotes the initial state. 
\begin{subequations}
    \begin{align}
        p(x_{t+1}|\xi_t) &= \delta (x_{t+1} - f_t(\xi_t))\\
        p(x_0) &= \delta(x_0 - \mathfrak{x}_0)
    \end{align}
\end{subequations}

Second, we formally characterize the behaviour of any controller tasked with administering the actions by the probabilistic policies, $\pi_t(u_t|x_t)$. These policies determine the probability of taking an action, $u_t$, at time $t$ for a given state $x_t$. As a result, we can represent the probabilistic closed-loop system dynamics with the following joint density function
\begin{equation}\label{eq:factorization_closedloop1}
p(\xi_{0:T};\pi_{0:T-1}) = p(x_0)\prod\nolimits_{t = 0}^{T-1} p(x_{t+1}|\xi_t)\pi_t(u_t|x_t)
\end{equation}
which evaluates the probability of any state-action trajectory, $\xi_{0:T}$, for a given policy sequence, $\pi_{0:T-1}$. Given that the probabilistic policies are at the heart of the approach, we emphasize notationally, that the closed-loop state-action trajectory density is parametrized by the policy sequence, $\pi_{0:T-1}$.


Thirdly, we introduce an exponential cost transformation, specifically, we substitute $\exp(-\gamma C(\xi_{0:T}))$ where $\gamma \in \mathbb{R}^+$. These formal changes allow us to rewrite (\ref{eq:general_cost_functional}) as a stochastic optimal control problem with exponential cost measure or so-called Risk-Sensitive Optimal Control problem (RSOC).
\begin{equation} \label{eq:risk_sensitivity_objective}
	\min_{\pi_{0:T-1}} \underbrace{ - \log 
		\mathbb{E}_{p(\xi_{0:T}| \pi_{0:T-1})} \left[ \exp \left(-\gamma C\left( \xi_{0:T} \right) \right) \right]  }_{ \textstyle\triangleq J(\pi_{0:T-1};\gamma)}
\end{equation}

{
\begin{lemma}
	\label{lemma:DOC=RSOC}
	Problem (\ref{eq:general_min_cost_functional}) and (\ref{eq:risk_sensitivity_objective}) are equivalent for deterministic dynamics {regardless of the value of $\gamma\in\mathbb{R}^+_*$. }
\end{lemma}
\begin{proof}
	This is a direct result of the deterministic dynamics which allows to drop the expectation. The equivalence of the problems is then trivial. 
\end{proof}

}

The proposed cost transformation has two important consequences. First, it will allow us to treat the problem as if it were a Maximum Likelihood Estimation (see later). Second, it also introduces the notion of risk into the framework. Optimal control problems of the form (\ref{eq:risk_sensitivity_objective}) are said to be risk seeking for positive $\gamma$. Indeed one verifies that the parameter determines the effective weight that is attributed to the tail of the variable cost's distribution. Hence the larger $\gamma$ is, the more optimistic the controller is about the outcome.

\begin{remark}
	In what follows the variability of the cost will be a consequence of the deliberate variability of the policies. Hence it is anticipated that parameter $\gamma$ can be used to control the level of the injected variation into the nonlinear optimization problem and thus determine the level of exploratory behaviour that any probabilistic policies might excite. 
\end{remark}

	In this work, we seek to extend the DDP algorithm by exploiting the uncertainty on the policy to better probe the optimization landscape and consequently realize an improved convergence rate. However, as reflected by Lemma \ref{lemma:DOC=RSOC}, for now, there is no incentive to actually produce any probabilistic policy. Such will be realised in the next sections.


\begin{figure}[t] 
	\centering 
	\resizebox{1\columnwidth}{!}{
		\begin{tikzpicture}
			\tikzstyle{main}=[circle, minimum size = 8mm, thick, draw =black!80, node distance = 6mm]
			\tikzstyle{connect} = [-latex,thick]
			\tikzstyle{floating} = [thick,dotted]
			\tikzstyle{measure} = [-latex,thick,bend right=45]
			
			\node[main] (X0) {$x_0$};
			\node[main] (X1) [right=of X0] {$x_1$};
			\node[main] (Xt) [right=of X1] {$x_t$};
			\node[main] (Xtt) [right=of Xt,label=center:$x_{t+1}$] {};
			\node[main] (XT) [right=of Xtt] {$x_{T}$};
			\node[main,draw=white] (XTT) [right=of XT] {};

			\node[main] (U0) [above left=of X1] {$u_0$};
			\node[main] (U1) [right=of U0] {$u_1$};
			\node[main] (Ut) [right=of U1] {$u_t$};
			\node[main] (Utt) [right=of Ut,label=center:$u_{t+1}$] {};
			
			\node[main,fill=gray!25] (Z0) [above left=of U1] {$\mathcal{O}_0$};
			\node[main,fill=gray!25] (Z1) [right=of Z0] {$\mathcal{O}_1$};
			\node[main,fill=gray!25] (Zt) [right=of Z1] {$\mathcal{O}_t$};
			\node[main,fill=gray!25] (Ztt) [right=of Zt,label=center:$\mathcal{O}_{t+1}$] {};
			\node[main,fill=gray!25] (ZT) [right=of Ztt] {$\mathcal{O}_{T}$};
			
			\node[main,draw=white] (n1) [above left=of U0] {};
			
			\path (X0) edge [connect] (X1);
			(X1) edge [floating] (Xt);
			(Xt) edge [connect] (Xtt);
			(Xtt) edge [floating] (XT);
			
			\path (X0) edge [measure] (Z0);
			\path (X1) edge [measure] (Z1);
			\path (Xt) edge [measure] (Zt);
			\path (Xtt) edge [measure] (Ztt);

            \path (XT) edge [connect] (ZT);
			
			\path (X0) edge [connect] (U0);
			\path (X1) edge [connect] (U1);
			\path (Xt) edge [connect] (Ut);
			\path (Xtt) edge [connect] (Utt);
			
			\path (U0) edge [connect] (X1);
			\path (U1) edge [floating] (Xt);
			\path (Ut) edge [connect] (Xtt);
			\path (Utt) edge [floating] (XT);
			
			\path (U0) edge [connect] (Z0);
			\path (U1) edge [connect] (Z1);
			\path (Ut) edge [connect] (Zt);
			\path (Utt) edge [connect] (Ztt);
		\end{tikzpicture}
	}
	\caption{\small{The probabilistic model used to represent the optimal control problem. White-shaded variables are latent or hidden.}}\label{fig:PGM}
\end{figure}

\subsection{Equivalent Maximum Likelihood Estimation problem}\label{sec:equivalent-probabilistic-inference-problem}

The key to treating (\ref{eq:risk_sensitivity_objective}) as a probabilistic inference problem is to internalize the notion of reward into the probabilistic model that was used to represent the dynamics. To this end, we extend the Markov model to a hidden Markov model by introducing a sequence of dummy random variables, $\mathcal{O}_{0:T}$. We treat the dummy variables as we would any other observation in the context of such a model. The state-action trajectory variables are now regarded as latent variables. Fig. \ref{fig:PGM} visualizes the corresponding probabilistic model.



%
%
%


To establish the connection with optimal control the following observation model is proposed. When $\mathcal{O}_t = 1$ it is implied that time step $t$ is optimal, any other value indicates that it is not. To exactly reproduce the problem (\ref{eq:risk_sensitivity_objective}), we define that the probability of being optimal at time $t$ is proportional to an exponential utility transform of the immediate cost, $c_t(\xi_t)$. 
\begin{align}\label{eq:exp_cost}
p(\mathcal{O}_t=1|\xi_t;\gamma) = \eta(\gamma) \exp(-\gamma c_t (\xi_t ))
\end{align}
{Here $\eta(\gamma)$ is a normalization factor so that $\int p(\mathcal{O}_t|\xi_t;\gamma)\text{d}\mathcal{O}_t$ is equal to $1$. The ensuing developments prevent the normalization factor from depending on $\xi_t$. An obvious choice is to choose $p(\mathcal{O}_t|\xi_t;\gamma)$ from the family of exponential functions so that $\eta(\gamma)=\gamma$. This is equivalent to assuming that $\int p(\mathcal{O}_t|\xi_t;\gamma)\text{d}\mathcal{O}_t = \eta(\gamma)\int_0^\infty \exp(-\gamma c) \text{d}c$ \cite{Peters2007}. This choice will have a desirable effect later on.} Finally, note that we will write $\mathcal{O}_t$ when we mean $\mathcal{O}_t=1$ henceforth.

Corresponding these adjustments, one verifies that the joint density of the full probabilistic model is now given by
\begin{multline}\label{eq:factorization}
p(\xi_{0:T},\mathcal{O}_{0:T};\pi_{0:T-1},\gamma) = p(x_0)p(\mathcal{O}_T|x_T;\gamma)\\\times \prod\nolimits_{t = 0}^{T-1} p(\mathcal{O}_t|\xi_t;\gamma)p(x_{t+1}|\xi_t)\pi_t(u_t|x_t)
\end{multline}
where $p(\mathcal{O}_T|x_T;\gamma)$ is associated with the terminal cost and is defined similarly as $p(\mathcal{O}_t|\xi_t;\gamma)$.


Now the connection with (\ref{eq:risk_sensitivity_objective}) can be established \cite{rawlik2013probabilistic,Noorani2022,lefebvre2023}.

\begin{theorem}\label{th:RSOC=MLE}
The minimization of the risk-sensitive objective function (\ref{eq:risk_sensitivity_objective}) is equivalent to maximizing the log-likelihood of the optimality variables, $\mathcal{O}_{0:T}$, that is
\begin{equation}\label{eq:Th_max_likelihood}
	\argmin_{\pi_{0:T-1}} \; J(\pi_{0:T-1};\gamma) = \argmax_{\pi_{0:T-1}} \; \log p(\mathcal{O}_{0:T};\pi_{0:T-1},\gamma)
\end{equation}
\end{theorem}

\begin{proof}
First, we can rewrite (\ref{eq:factorization}) in terms of (\ref{eq:factorization_closedloop1}). Then we substitute the definition (\ref{eq:exp_cost}) for the measurement models. This yields
\begin{equation*}\label{eq:factorization2}
\begin{multlined}
p(\xi_{0:T}, \mathcal{O}_{0:T};\pi_{0:T-1},\gamma)  \\=
\gamma^{T+1} p(\xi_{0:T};\pi_{0:T-1}) \exp\left(-\gamma C\left(\xi_{0:T} \right)\right)
\end{multlined}
\end{equation*}
Marginalizing for $\xi_{0:T}$ and taking the logarithm retrieves the log-likelihood of the optimality variables. By construction, this expression turns out to be equivalent to the negative risk-sensitive objective function in (\ref{eq:risk_sensitivity_objective}) plus some constant. Therefore the extremum points of both objectives are equivalent. 
\begin{equation*} \label{eq:max_likelihood}
    \begin{multlined}
    \log p(\mathcal{O}_{0:T};\pi_{0:T-1},\gamma)  \\
    \begin{aligned}
        &= \log \int p(\xi_{0:T}, \mathcal{O}_{0:T};\pi_{0:T-1},\gamma) \text{d}\xi_{0:T} \\
        &= \log \int p(\xi_{0:T};\pi_{0:T-1}) \exp\left(-\gamma C\left(\xi_{0:T} \right)\right) \text{d}\xi_{0:T} +\log\gamma^{T+1}
    \end{aligned}
    \end{multlined}
    \vspace*{-12pt}
\end{equation*}
\end{proof}

From the theorem above it directly follows that we can treat the problem on the right-hand side of (\ref{eq:Th_max_likelihood}) instead of (\ref{eq:risk_sensitivity_objective}). MLE problems of this form are usually addressed directly using the EM algorithm. 

\subsection{Expectation-Maximization}\label{sec:expectation-maximization}

The EM algorithm is a well-established technique to transform hard MLE problems into a sequence of easier subproblems where the result of every subproblem serves as the input to the next subproblem. Instead of maximizing the MLE directly, in each subproblem, a lower bound on the log-likelihood of the observations is maximized. This lower bound is referred to as the evidence lower bound (ELBO). 

To practice the EM algorithm we need to define the ELBO for the log-likelihood of the observations. The following lemma states the ELBO corresponding with the right-hand side of the equation (\ref{eq:Th_max_likelihood}) \cite{bishop2006pattern}. 


\begin{lemma}
    The following is true for any density, $q(\xi_{0:T})$
    \begin{equation} 
        \begin{multlined}
            \log  p(\mathcal{O}_{0:T};\pi_{0:T-1},\gamma)  \\
            \label{eq:ELBO_KL}
                = \underbrace{ \mathbb{E}_{q(\xi_{0:T})} \left[ \log \frac{ p(\xi_{0:T}, \mathcal{O}_{0:T};\pi_{0:T-1},\gamma)}{q(\xi_{0:T})} \right]}_{\textstyle\triangleq \mathcal{E}(q,\pi_{0:T-1},\gamma)} \\
                +\mathbb{KL}\left[ q(\xi_{0:T}) \parallel p(\xi_{0:T}|\mathcal{O}_{0:T};\pi_{0:T-1},\gamma) \right]
        \end{multlined}
    \end{equation}
Here $\mathbb{KL}$ denotes the Kullback–Leibler (KL) divergence (or relative entropy) and $\mathcal{E}$ denotes the ELBO. Due to the positive definiteness of the KL divergence, it follows that
\begin{equation}
    \begin{aligned} \label{eq:ELBO}
        \log  p(\mathcal{O}_{0:T};\pi_{0:T-1},\gamma)  \ge
        \mathcal{E}(q,\pi_{0:T-1},\gamma)
    \end{aligned} 
\end{equation}
\end{lemma}
\begin{proof}
	First note that the KL-divergence is defined as
	\begin{equation*}
		\mathbb{KL}[q||p] = \int q(x)\log\frac{p(x)}{q(x)}\text{d}x
	\end{equation*}
    Then, by expanding the right-hand side of (\ref{eq:ELBO_KL}), and relying on the sum of logarithms property and Bayes’ rule, one obtains
    \begin{multline*}
            \int q(\xi_{0:T})
            \log \left(\frac{p(\xi_{0:T}, \mathcal{O}_{0:T};\pi_{0:T-1},\gamma)}{q(\xi_{0:T})}\right. 
            \\
            \left.\times\frac{q(\xi_{0:T})}{p(\xi_{0:T}| \mathcal{O}_{0:T};\pi_{0:T-1},\gamma)} \right)\text{d}\xi_{0:T} \\ =\log p(\mathcal{O}_{0:T};\pi_{0:T-1},\gamma) 
    \end{multline*} 
\end{proof}

Indeed the ELBO establishes a lower bound on the log-likelihood of the observations. Then in every main iteration of the EM algorithm, we need to execute the following two steps commonly referred to as the E- and M-step.
\begin{itemize}
\item \textit{Expectation step (E-step):} The E-step aims to minimize the gap between the ELBO and the log-likelihood of the observations by choosing the optimal inference density, $q(\xi_{0:T})$. The ELBO (\ref{eq:ELBO}) is highest when the KL divergence term in (\ref{eq:ELBO_KL}) is lowest, more specifically equal to zero. This happens when the inference density, $q(\xi_{0:T})$, is equivalent to the posterior density, $p(\xi_{0:T}|\mathcal{O}_{0:T}; \pi_{0:T-1},\gamma)$. Note that the parameters of the inference density are fixed. To avoid confusion with the optimization parameters, we substitute $\rho_{0:T-1}$ and $\lambda$ for $\pi_{0:T-1}$ and $\gamma$, respectively.
\begin{equation}\label{eq:Estep}
    q^*(\xi_{0:T}) = p(\xi_{0:T}|\mathcal{O}_{0:T}; \rho_{0:T-1},\lambda)
\end{equation}
\item \textit{Maximization step (M-Step):}
The M-step then maximizes the ELBO for given optimal inference density, $q^*$, and with respect to the parameters $\pi_{0:T-1}$ and $\gamma$. This problem can be solved explicitly as we will show below.
\begin{equation}\label{eq:Mstep}
        \{\pi^*_{0:T-1},\gamma^*\} = \argmax_{\pi_{0:T-1},\gamma} \; \mathcal{E}(q^*,\pi_{0:T-1},\gamma)
\end{equation}
\end{itemize}

After the M-step is completed, the algorithm repeats the steps, substituting the optimal policies, $\pi_{0:T-1}^*$, and parameter, $\gamma^*$, for the prior policies, $\rho_{0:T-1}$, and parameters, $\lambda$, in the subsequent E-step. The EM algorithm thus establishes a fixed-point iteration producing a sequence of probabilistic policies. 

\begin{lemma}\label{lem:EM}
    The EM algorithm is guaranteed to converge to a local optimum of the corresponding MLE problem.
\end{lemma}
\begin{proof}
		Refer to Chapter 6.5.3 \cite{murphy2023}. 
\end{proof}
{\begin{remark}
	Although maximizing for $\gamma$ in the M-step is not strictly required to apply the EM algorithm to solve the problem (\ref{eq:risk_sensitivity_objective}), there is no formal objection not to. The equivalence between (\ref{eq:general_min_cost_functional}) and (\ref{eq:risk_sensitivity_objective}) is unaffected by the value of $\gamma$ as was established by Lemma \ref{lemma:DOC=RSOC}. 
    Note that this is only the case for deterministic systems. For stochastic systems, this would result into the solution of another RSOC problem.
\end{remark}}
The temperature update was first considered in a stochastic setting by \cite{watson2020stochastic,watson2021b} and dates back to \cite{Peters2007} for the given cost transformation.

Let us now further evaluate the M-step. The following theorem expresses the explicit solution of (\ref{eq:Mstep}). {
}

\begin{theorem} \label{th:EM_solution}
    The objective in (\ref{eq:Mstep}) attains an optimum for
    \begin{subequations}
        \begin{align}
\label{eq:Opt_Policy}
            \pi^*_{t} &= \frac{p(\xi_{t}|\mathcal{O}_{0:T}; \rho_{0:T-1},\lambda)}{p(x_{t}|\mathcal{O}_{0:T}; \rho_{0:T-1},\lambda)} = p(u_{t}|x_{t}, \mathcal{O}_{0:T}; \rho_{0:T-1},\lambda)\\
            \gamma^* &= \frac{T+1}{\mathbb{E}_{p(\xi_{0:T}|\mathcal{O}_{0:T}; \rho_{0:T-1},\lambda)}\left[ C(\xi_{0:T}) \right]}
            \label{eq:Opt_Parameter}
        \end{align}
    \end{subequations}
    Here $p(u_{t}|x_{t}, \mathcal{O}_{0:T}; \rho_{0:T-1},\lambda)$ represents the conditional probability of the action $u_t$ conditioned on the state, $x_t$, and the optimality variables, $\mathcal{O}_{0:T}$, parametrised by the former policy sequence, $\rho_{0:T-1}$, and risk-sensitivity parameter, $\lambda$.
\end{theorem}
\begin{proof}
    By substituting the factorisation mentioned in (\ref{eq:factorization}) into the ELBO and using the logarithm of a product property, we can rewrite the optimization problem (\ref{eq:Mstep}) as follows:
    \begin{equation*}
    	\begin{multlined}
    		    	\pi^*_{0:T-1} \\
    		\begin{aligned}
        &= \argmax_{\pi_{0:T-1}} \; \mathcal{E}(q^*,\pi_{0:T-1},\gamma) \nonumber\\ 
        &= \argmax_{\pi_{0:T-1}} \mathbb{E}_{p(\xi_{0:T}|\mathcal{O}_{0:T}; \rho_{0:T-1},\lambda)} \left[\sum\nolimits_{t=0}^{T-1} \log  \pi_t(u_t|x_t) \right] \nonumber\\
        &= \argmax_{\pi_{0:T-1}} \sum\nolimits_{t=0}^{T-1} \mathbb{E}_{p(\xi_{t}|\mathcal{O}_{0:T}; \rho_{0:T-1},\lambda)} \left[ \log  \pi_t(u_t|x_t) \right] 
    \end{aligned}
    	\end{multlined}
    \end{equation*}
    Then we can apply the calculus of variations on the following dynamic optimization problem to find the optimal policy 
    \begin{equation}
        \begin{aligned}
            \pi^*_{t} = &\argmax_{\pi_{t}} \int p(\xi_t|\mathcal{O}_{0:T}; \rho_{0:T-1},\lambda)  \log  \pi_t(u_t|x_t) \text{d}\xi_{t}
            \nonumber \\
            &\text{ s.t.} \int \pi_t(u_t|x_t) \text{d}u_t=1 \nonumber
        \end{aligned} 
    \end{equation}
    
    One verifies
    \begin{equation*}
        \pi^*_{t} = \frac{p(\xi_{t}|\mathcal{O}_{0:T}; \rho_{0:T-1},\lambda)}{\int p(\xi_{t}|\mathcal{O}_{0:T}; \rho_{0:T-1},\lambda)\text{d}u_t}
    \end{equation*}
    Indeed the denominator marginalizes over the action, $u_t$, resulting into $p(x_{t}|\mathcal{O}_{0:T}; \rho_{0:T-1},\lambda)$. The second equality follows from Bayes' rule. 
    
    Likewise, we can maximize the ELBO with respect to $\gamma$.
    \begin{equation*}
        \begin{multlined}
        	\gamma^*\\ 
        	\begin{aligned}
             &= \argmax_{\gamma} \mathcal{E}(q^*,\pi_{0:T-1},\gamma) \nonumber\\ 
            &= \argmax_{\gamma} \mathbb{E}_{p(\xi_{0:T}|\mathcal{O}_{0:T}; \rho_{0:T-1},\lambda)} \left[\sum\nolimits_{t=0}^{T} \log  p(\mathcal{O}_t|\xi_t;\gamma) \right] \nonumber\\
            &= \argmax_{\gamma} \; (T+1) \log \gamma - \gamma \mathbb{E}_{p(\xi_{t}|\mathcal{O}_{0:T};\rho_{0:T-1},\lambda)}\left[ C(\xi_{0:T}) \right]
        \end{aligned}
        \end{multlined}
    \end{equation*}
    Taking the derivative of the above expression with respect to $\gamma$ and setting it to zero retrieves the required expression. 
\end{proof}

\begin{remark}
    Maximizing the ELBO with respect to $\gamma$ is equivalent to minimising the risk-sensitive objective function (\ref{eq:risk_sensitivity_objective}) with respect to $\gamma$ plus the regularization term $-(T+1)\log\gamma$.
\end{remark}    

We emphasize that if it were not for the normalization factor in (\ref{eq:exp_cost}), the trivial solution to minimize the risk-sensitive objective function (\ref{eq:risk_sensitivity_objective}) with respect to $\gamma$ would have been zero. Minimizing the objective function (\ref{eq:risk_sensitivity_objective}) with respect to $\gamma$ can help to control how risk-seeking the corresponding policies are.
The regularization term $-(T+1)\log\gamma$ prevents the optimal value of the risk-seeking parameter from reaching zero. As a result, there exists an inverse relationship between the total expected cost and the value of the optimal parameter, $\gamma^*$. Hence, if the prior policy sequence, $\rho_{0:T-1}$, reduces the expectation of the total cost, then the risk parameter will be increased in the next iteration, thus increasing the optimism of the optimal policies. This property is highly desirable as it counteracts the imminent collapse of the policies. Following Theorem \ref{th:RSOC=MLE} and Lemma \ref{lem:EM}, the policies are bound to slowly converge to deterministic policies as the EM iterations increase. In some situations, this might result in a premature convergence of the policies. Actively stimulating exploration as the EM algorithm converges can address this problem. 



\subsection{Evaluation of the optimal policies}\label{sec:iterative-solution-methods}

As was shown in the previous section, treatment of the MLE (\ref{eq:Th_max_likelihood}) using the EM algorithm generates a sequence of probabilistic policies that are guaranteed to converge to the deterministic optimal control defined in (\ref{eq:general_min_cost_functional}). It remains to show how we can evaluate the optimal policies in (\ref{eq:Opt_Policy}) efficiently. Here we discuss two possible evaluation strategies. The first approach exhibits a close resemblance with the principle of dynamic programming. The second approach leans on well-established techniques to practice inference on probabilistic models such as the one visualized in Fig. \ref{fig:PGM}.

\subsubsection{Probabilistic Dynamic Programming} \label{sec:PDP}
The following theorem establishes a dynamic programming approach to evaluate the optimal probabilistic policy (\ref{eq:Opt_Policy}), e.g. \cite{levine2018,lefebvre2023}.
\begin{theorem} \label{th:IterativeProbabilisticDP}
    The optimal policy sequence defined in (\ref{eq:Opt_Policy}) can be represented as
    \begin{equation} \label{eq:PDP_Policy}
    	\pi_t^*(u_t|x_t) = {\rho}_{t}(u_t|x_t) \frac{\exp(-\lambda {Q}^*_t(\xi_t))}{\exp(-\lambda {V}^*_t(x_t))}
    \end{equation}
    in which the probabilistic state-action value function, $Q^*_t$, and, optimal cost-to-go, $V^*_t$, are defined as
    \begin{subequations} \label{eq:PDP_ValueFunctions_Definitions}
        \begin{align} 
            {Q}^*_t(\xi_t) &\triangleq -\tfrac{1}{\lambda}\log p(\mathcal{O}_{t:T}|\xi_{t}; \rho_{t+1:T-1},\lambda) \\
            {V}^*_t({x}_t) &\triangleq -\tfrac{1}{\lambda}\log p(\mathcal{O}_{t:T}|x_{t}; \rho_{t:T-1},\lambda)
        \end{align}
    \end{subequations}
  {The functions satisfy the following backward recursion
    \begin{subequations} \label{eq:PDP_ValueFunctions}	
        \begin{align} 
        	 \label{eq:PDP_ValueFunctions_Q}
            Q^*_{t}(\xi_t) &= \bar{c}_t({\xi}_t)-\tfrac{1}{\lambda}\log \expect{p({x}_{t+1}|{\xi}_t)}\left[\exp(-\lambda V^*_{t+1}({x}_{t+1}))\right] \\
            V^*_{t}({x}_t) &= -\tfrac{1}{\lambda} \log \expect{\rho_t(u_t|x_t)}\left[\exp(-\lambda {Q}^*_t(\xi_t))\right]
             \label{eq:PDP_ValueFunctions_V}
        \end{align}
    \end{subequations}
    where $\bar{c}_t = c_t - \frac{1}{\lambda} \log \lambda$. The recursion is initialised with $V^*_T = \bar{c}_T(x_T)$.}
\end{theorem}
\begin{proof}
    First note that, as illustrated by the model in Fig. \ref{fig:PGM}, the action $u_t$ and $\mathcal{O}_{0:t-1}$ are conditionally independent given $x_t$. It follows that $p(u_{t}|x_{t}, \mathcal{O}_{0:T}; \rho_{0:T-1},\lambda)$ simplifies to $p(u_{t}|x_{t}, \mathcal{O}_{t:T}; \rho_{t:T-1},\lambda)$. Then using Bayes' rule, we can rewrite the optimal policy (\ref{eq:Opt_Policy}) to establish (\ref{eq:PDP_Policy}).  
    \begin{equation*}
    	\begin{multlined}
    		\pi_t^*(u_t|x_t) \\
    		\begin{aligned}
            &= p(u_{t}|x_{t}, \mathcal{O}_{t:T}; \rho_{t:T-1},\lambda) \\
            &= \frac{p(\xi_{t}|\mathcal{O}_{t:T}; \rho_{t:T-1})}{p(x_{t}|\mathcal{O}_{t:T}; \rho_{t:T-1},\lambda)} \\ 
            &= \frac{{\rho}_{t}(u_t|x_t) p(\mathcal{O}_{t:T}|\xi_{t}; \rho_{t+1:T-1},\lambda)}{p(\mathcal{O}_{t:T}|x_{t}; \rho_{t:T-1},\lambda)}
        \end{aligned}
    	\end{multlined}
    \end{equation*}
     
     Second, to find the relation between, $Q^*_t$, and, $V^*_t$, we can factorise $p(u_t,\mathcal{O}_{t:T}|x_t;\rho_{t:T-1},\lambda)$ 
     \begin{equation*}
\begin{multlined}
	     	p(u_t,\mathcal{O}_{t:T}|x_t;\rho_{t:T-1},\lambda) = \rho_t(u_t|x_t) p(\mathcal{O}_{t:T}|\xi_t;\rho_{t+1:T-1},\lambda)
\end{multlined}
     \end{equation*}
    Then we marginalize out $u_t$ to obtain
    \begin{equation*}
       \begin{multlined}
       	p(\mathcal{O}_{t:T}|x_{t}; \rho_{t:T-1},\lambda) \\
       	 \begin{aligned}
            &= \int \rho_t(u_t|x_t) p(\mathcal{O}_{t:T}|\xi_{t}; \rho_{t+1:T-1},\lambda) \text{d}u_t \\
            & = \expect{\rho_t(u_t|x_t)}\left[p(\mathcal{O}_{t:T}|\xi_{t}; \rho_{t+1:T-1},\lambda)\right]
        \end{aligned}
       \end{multlined}
    \end{equation*}
     
     Thirdly, to find the recursive expression for $V_t^*$, we observe that
     \begin{equation*}
\begin{multlined}
	     	p(\mathcal{O}_{t:T}|\xi_{t}; \rho_{t+1:T-1},\lambda) \\= p(\mathcal{O}_{t}|\xi_t;\lambda) p(\mathcal{O}_{t+1:T}|\xi_{t}; \rho_{t+1:T-1},\lambda)
\end{multlined}
     \end{equation*}
     Now we factorise $p(x_{t+1},\mathcal{O}_{t+1:T}|\xi_t;\rho_{t+1:T-1},\lambda)$
     \begin{equation*}
\begin{multlined}
	     	p(x_{t+1},\mathcal{O}_{t+1:T}|\xi_t;\rho_{t+1:T-1},\lambda) \\ = p(x_{t+1}|\xi_t)p(\mathcal{O}_{t+1:T}|x_{t+1};\rho_{t+1:T-1},\lambda)
\end{multlined}
     \end{equation*}
     Combination with the previous result and marginalization of $x_{t+1}$ shows that
    \begin{equation*}
    	\begin{multlined}
    		p(\mathcal{O}_{t:T}|\xi_{t}; \rho_{t+1:T-1},\lambda) \\
        \begin{aligned}
            &= p(\mathcal{O}_{t}|\xi_t;\lambda) \expect{p(x_{t+1}|\xi_t)} \left[ p(\mathcal{O}_{t+1:T}|x_{t+1}; \rho_{t+1:T-1},\lambda) \right]
        \end{aligned}
    	\end{multlined}
    \end{equation*}
    
    Finally, taking the logarithm of the last second and third results and multiplying by $-\frac{1}{\gamma}$, we recover the desired expressions.
\end{proof}

\begin{remark}\label{remark:PDP_Q}
    For a deterministic dynamic system, the recursive expression for the probabilistic state-action value function in (\ref{eq:PDP_ValueFunctions}) collapses on the state-action value function defined in (\ref{eq:Q_function}) {\color{black}up to a constant $\frac{1}{\gamma} \log \gamma$.}
\end{remark}


\subsubsection{Bayesian Smoothing Control} \label{sec:BSC}

Evaluation of the posterior probability, $p(\xi_{0:T}|\mathcal{O}_{0:T}; \rho_{0:T-1},\lambda)$, and in particular $p(\xi_t|\mathcal{O}_{0:T}; \rho_{0:T-1},\lambda)$, also happens to coincide with the Bayesian smoother. The Bayesian smoother determines the probability of the variable $\xi_{t}$ conditioned on the variables $\mathcal{O}_{0:T}$. If $p(\xi_t|\mathcal{O}_{0:T}; \rho_{0:T-1},\lambda)$ were available, the optimal policies could be evaluated simply by practising the definition in (\ref{eq:Opt_Policy}). The Bayesian smoother is a well-established inference problem and can be evaluated using the Bayesian smoothing equations. The following theorem is adopted from \cite{sarkka2023,murphy2023} and states the Bayesian smoothing equations adapted to the generalised probabilistic model in Fig. \ref{fig:PGM}. 

\begin{theorem}\label{th:BayesianInferenceControl}
    The backward recursive equation for computing the {smoothing densities} is given by 
    \begin{equation}
    \begin{multlined} \label{eq:BIC_URTSS1}
            p(\xi_{t}|\mathcal{O}_{0:T}; \rho_{0:T-1},\lambda) = p(\xi_t|\mathcal{O}_{0:t}; \rho_{0:t},\lambda) \\
            \times \int \frac{p(x_{t+1}|\xi_t)p(x_{t+1}|\mathcal{O}_{0:T}; \rho_{0:T-1},\lambda)}{p(x_{t+1}|\mathcal{O}_{0:t}; \rho_{0:t},\lambda)}\text{d}x_{t+1}
    \end{multlined}   
    \end{equation}
    where
    \begin{equation} \label{eq:BIC_URTSS2}
            p(x_{t+1}|\mathcal{O}_{0:t}; \rho_{0:t},\lambda) = \int p(x_{t+1}|\xi_{t}) p(\xi_t|\mathcal{O}_{0:t}; \rho_{0:t},\lambda) \text{d} \xi_t 
    \end{equation}
    
    The backward smoothing correction equations presume knowledge of the probability $p(\xi_t|\mathcal{O}_{0:t}; \rho_{0:t})$. This density is known as the {filtering density} and can be computed by the following forward recursive equation or so-called forward filtering equation
    \begin{equation}
    \begin{aligned}    \label{eq:BIC_UKF_correct}
            p(\xi_t|\mathcal{O}_{0:t}; \rho_{0:t},\lambda) = \frac{p(\mathcal{O}_{t}|\xi_{t};\lambda)p(\xi_t|\mathcal{O}_{0:t-1}; \rho_{0:t},\lambda)}{\int p(\mathcal{O}_{t}|\xi_{t};\lambda)p(\xi_t|\mathcal{O}_{0:t-1}; \rho_{0:t},\lambda) \text{d}\xi_t}
    \end{aligned}
    \end{equation}
    where 
    \begin{equation}
    \begin{multlined}  \label{eq:BIC_UKF_pred}
            p(\xi_t|\mathcal{O}_{0:t-1}; \rho_{0:t},\lambda) \\ =
            \int p(\xi_t|\xi_{t-1};\rho_{t}) p(\xi_{t-1}|\mathcal{O}_{0:t-1}; \rho_{0:t-1},\lambda)\text{d}\xi_{t-1}
    \end{multlined}  
    \end{equation}
    The density $p(\xi_t|\xi_{t-1};\rho_{t})$ is defined as the state-action transition probability of the dynamic system governed by the prior controller $\rho_t$. {Equations (\ref{eq:BIC_UKF_correct}) and (\ref{eq:BIC_UKF_pred}) are known as the update and prediction steps, respectively.}
    
    Finally, note that the forward recursion is initialised with
    \begin{equation}
    	p(\xi_0|\mathcal{O}_0;\rho_0,\lambda) \triangleq p(\xi_0;\rho_0) = p(x_0) \rho_0(u_0|x_0)
    \end{equation}
\end{theorem}
\begin{proof}
        Refer to Chapter 8.1 of \cite{sarkka2023}. 
\end{proof}

\section{Numerical implementation} \label{sec:algorithm}

	At this point, we may want to summarize what has been discussed so far. In the previous section, the nonlinear discrete-time deterministic optimal control problem (\ref{eq:general_min_cost_functional}) was formulated as an equivalent MLE problem (\ref{eq:Th_max_likelihood}). This allowed to treat the problem (\ref{eq:general_min_cost_functional}) using the EM algorithm. Treatment of the MLE using the EM algorithm had two important consequences. The original problem reduces to a sequence of recursive subproblems. As a result, the method is iterative by design. Second, the extremum of every subproblem is attained for a probabilistic policy. It is here hypothesized that the uncertainty of the policy will scale with the degree of exploration that is required at that stage of the solution. 
	
	However, so far section \ref{sec:inference} has only established theory in the same sense that the principle of dynamic programming (recall section \ref{sec:dynamic-programming}) establishes the theory for the DDP algorithm in section \ref{sec:differential-dynamic-programming}. In this section, we discuss how the EM treatment of the MLE problem can be numerically applied to practical problems. Since there are two evaluation approaches for the optimal policies in (\ref{eq:Opt_Policy}), see Theorems  \ref{th:IterativeProbabilisticDP} and \ref{th:BayesianInferenceControl}, it is possible to derive two distinct algorithms. Both algorithms turn out to be gradient-free alternatives to DDP that can be shown to converge faster than the traditional DDP algorithm, empirically. 


Provided the probabilistic nature of the entire framework, the numerical algorithms rely on a considerable amount of Uncertainty Propagation (UP). State-of-the-art methods for UP are Gaussian moment matching in combination with sigma-point methods \cite{sarkka2023}. Gaussian moment matching is a method used to approximate a complex or non-Gaussian distribution with a Gaussian distribution by ensuring that the first and second moments (mean and covariance) of the Gaussian match those of the original distribution. The sigma-point method is a popular technique to approximately evaluate the expected value of a nonlinear map under Gaussian measures.

\begin{definition} [Gaussian moment matching \cite{sarkka2023}] 
    The Gaussian moment matching approximation to the joint probability of the random variables $\alpha\in\mathbb{R}^{n_\alpha}$ and the transformed random variable $\beta=g(\alpha)$, where $g:\mathbb{R}^{n_\alpha}\mapsto\mathbb{R}^{n_\beta}$ is some nonlinear map and $\alpha \sim \mathcal{N}(\alpha;\mu_\alpha, \Sigma_{\alpha\alpha})$, is defined as    
    \begin{equation} \label{eq:Moment_Matching1}
        \begin{aligned}
               \begin{bmatrix}
               	\alpha \\\beta
               \end{bmatrix} 
               \sim \mathcal{N} 
                \left(
                   \begin{bmatrix}
                   	\alpha \\\beta
                   \end{bmatrix};\begin{bmatrix}
                   \mu_\alpha \\ \mu_\beta
                   \end{bmatrix} ,\begin{bmatrix}
\Sigma_{\alpha\alpha}& \Sigma_{\alpha\beta} \\ 
* & \Sigma_{\beta\beta}
                   \end{bmatrix} 
                \right)
        \end{aligned}
    \end{equation}
    where
    \begin{subequations} \label{eq:Moment_Matching2}
        \begin{align}
            \mu_{\beta} &= \expect{\mathcal{N}(\alpha;\mu_{\alpha}, \Sigma_{\alpha \alpha})} \left[ g(\alpha)\right] \\
            \Sigma_{\beta \beta} &= \expect{\mathcal{N}(a;\mu_{\alpha}, \Sigma_{\alpha \alpha})} \left[ (g(\alpha)-\mu_{\beta})(g(\alpha)-\mu_{\beta})^{\top}\right]\\
            \Sigma_{\alpha \beta} &= \expect{\mathcal{N}(\alpha;\mu_{\alpha}, \Sigma_{\alpha \alpha})} \left[ (\alpha-\mu_{\alpha})(g(\alpha)-\mu_{\beta})^{\top}\right]
        \end{align}
    \end{subequations}
\end{definition}

\begin{definition} [Sigma-point method \cite{sarkka2023}] \label{recall:SP}
    The sigma-point method approximates the expectation of the nonlinear map, $g$, under the measure, $\mathcal{N}(\alpha;\mu_\alpha,\Sigma_{\alpha\alpha})$ with
    \begin{equation} \label{eq:Sigma_Points}
        \begin{aligned}
            \expect{\mathcal{N}(\alpha;\mu_{\alpha},\Sigma_{\alpha \alpha})}\left[ g(\alpha) \right] &= \int g(\alpha) \mathcal{N}(\alpha|\mu_{\alpha},\Sigma_{\alpha \alpha}) d\alpha\\
            &\approx \sum\nolimits_{n=1}^{N_\alpha} w^{\alpha}_n g(\sqrt{\Sigma}_{\alpha \alpha} \epsilon^{\alpha}_n + \mu_{\alpha})
        \end{aligned}
    \end{equation}
    where $\sqrt{\Sigma}_{\alpha \alpha}$ denotes the Cholesky factor or some other similar square root of the covariance matrix $\Sigma_{\alpha \alpha}$, $N_{\alpha}$ is the number of the required sigma points, $\epsilon^{\alpha}_n \in \mathbb{R}^{n_{\alpha}}$ refers to the $n_{\alpha}$-dimensional unit sigma points with $n\in\{1,...,N_{\alpha}\}$ and $w^{\alpha}_n$ are their associated scalar weights. The determination of the signature relies on the predetermined sigma-point method and can be done in several ways. There are two main strategies to generate the sigma points related to the Gauss–Hermite and spherical cubature integration methods \cite{sarkka2023}.
\end{definition}

\begin{definition} [Unscented transform\cite{sarkka2023}] \label{recall:UT}
    The unscented transform refers to the combination of the Gaussian moment matching and the sigma-point methods. Trivially one uses (\ref{eq:Sigma_Points}) to numerically approximate the integrals (\ref{eq:Moment_Matching2}). 
\end{definition}

Further, we also introduce the Fourier-Hermite series, which is closely related to the definitions above. The Fourier-Hermite series provides an alternative for the Taylor series to produce a second-order approximation of any nonlinear function that is approximately valid in a larger neighbourhood. 

\begin{definition} [Fourier-Hermite Series \cite{Sarmavuori2012}] \label{recall:FH}
	The second-order Fourier-Hermite series of the nonlinear map, $g$, is defined as 
	\begin{equation}
		g(\alpha) \approx \hat{g}(\alpha)=\frac{1}{2} \begin{bmatrix}
			1 \\ \delta \alpha
		\end{bmatrix}^\top \begin{bmatrix}
			2\hat{g}_0 & \hat{g}_{\alpha}^{\top}  \\
			* & \hat{g}_{\alpha \alpha}  \\
		\end{bmatrix}\begin{bmatrix}
			1 \\ \delta \alpha
		\end{bmatrix}
	\end{equation} 
	where 
	\begin{subequations}
		\begin{align}
			\hat{g}_0 &= \tilde{g}_0 - \tfrac{1}{2}\trace\{\Tilde{g}_{\alpha \alpha}\}\\
			\hat{g}_{\alpha} &= \sqrt{\Sigma}^{-1}_{\alpha\alpha}\tilde{g}_{\alpha}\\
			\hat{g}_{\alpha \alpha} &= \sqrt{\Sigma}_{\alpha\alpha}^{-1} \tilde{g}_{\alpha \alpha}\sqrt{\Sigma}_{\alpha\alpha}^{-1}
		\end{align}
	\end{subequations}
	and
	\begin{subequations}
		\begin{align}
			\tilde{g}_0 &= \expect{\mathcal{N}(\alpha;\mu_{\alpha}, \Sigma_{\alpha\alpha})}\left[g(\alpha)\right]\\
			\Tilde{g}_{\alpha} &=  \expect{\mathcal{N}(\alpha;\mu_{\alpha}, \Sigma_{\alpha \alpha})}\left[g(\alpha) H_1(\sqrt{\Sigma}^{-1}_{\alpha\alpha} \delta \alpha)\right]\\
			\tilde{g}_{\alpha \alpha} &= \expect{\mathcal{N}(\alpha;\mu_{\alpha}, \Sigma_{\alpha\alpha})}\left[g(\alpha) H_2(\sqrt{\Sigma}^{-1}_{\alpha\alpha} \delta \alpha)\right]
		\end{align}
	\end{subequations}
	where $\delta \alpha \triangleq \alpha - \mu_{\alpha}$, $H_0(\alpha) = 1$, $H_1(\alpha) = \alpha$ and $H_2(\alpha) = \alpha \alpha^T - \matrixstyle{I}$ are multivariate Hermite polynomials. The expectations are taken about an arbitrary Gaussian distribution $\mathcal{N}(\alpha;\mu_{\alpha}, \Sigma_{\alpha\alpha})$ such that $\mu_{\alpha}$ and $\Sigma_{\alpha}$ are indicative for some area of interest. The expectations can be evaluated using a sigma-point method as follows
	\begin{equation}
		\begin{aligned}
			\tilde{g}_0 &= \sum\nolimits_{n=1}^{N_\alpha} w^{\alpha}_n g\left(\mu_{\alpha}+\sqrt{\Sigma}_{\alpha\alpha}\epsilon^{\alpha}_n\right)\\
			\Tilde{g}_{\alpha} &= \sum\nolimits_{n=1}^{N_\alpha} w^{\alpha}_n g\left(\mu_{\alpha}+\sqrt{\Sigma}_{\alpha\alpha}\epsilon^{\alpha}_n\right) \epsilon^{\alpha}_n \\
			\tilde{g}_{\alpha \alpha} &= \sum\nolimits_{n=1}^{N_\alpha} w^{\alpha}_n g\left(\mu_{\alpha}+\sqrt{\Sigma}_{\alpha\alpha}\epsilon^{\alpha}_n\right) \left(\epsilon^{\alpha}_n\epsilon^{\alpha\top}_n -\matrixstyle{I}\right)
		\end{aligned}
	\end{equation}
\end{definition} 

In the next two sections, we derive two iterative numerical algorithms to address the problem in (\ref{eq:Th_max_likelihood}) and thus the original problem in (\ref{eq:general_min_cost_functional}). The main iterations are governed by the EM algorithm. Every iteration produces optimal probabilistic policies, $\pi_{0:T-1}^*$, which serve as the prior probabilistic policies, $\rho_{0:T-1}$, in the next iteration. Inspired by the DDP algorithm we propose to approximate the policies as follows
\begin{subequations} \label{eq:SPPDP_Cont}
	\begin{align}
		\label{eq:SPPDP_Cont_prior}
	 {\rho}_t &\approx \hat{\rho}_t= \mathcal{N}(\delta u_t;k_t+\matrixstyle{K}_t \delta x_t,\Sigma_t) \\
		 {\pi}_t^* &\approx\hat{\pi}_t^* = \mathcal{N}(\delta u_t;k_t^*+\matrixstyle{K}_t^* \delta x_t,\Sigma_t^*)
		\label{eq:SPPDP_Cont_posterior}
	\end{align}
\end{subequations}


\subsection{Sigma-Point Probabilistic Dynamic Programming}
The algorithm presented in this section derives directly from the policy evaluation procedure described in Theorem \ref{th:IterativeProbabilisticDP}.
The main idea is again similar to that of DDP. The goal is to obtain a quadratic approximation of the \textit{probabilistic} state-action value function around some nominal trajectory. The nominal trajectory is produced by forward closed-loop simulation of the prior probabilistic policies, $\rho_{0:T-1}$. In the present context, this produces a state-action trajectory probability rather than a unique trajectory. This invites to approximate the probabilistic state-action value function in  (\ref{eq:PDP_ValueFunctions_Q}) using a second-order Fourier-Hermite series expansion. In turn, the optimal cost-to-go can be evaluated by practising the definition (\ref{eq:PDP_ValueFunctions_V}). Finally, the optimal policy can be evaluated as in (\ref{eq:PDP_Policy}). The entire procedure is then repeated until convergence. 



\subsubsection{Forward pass}
As described above, we require access to an approximation of the nominal closed-loop trajectory density, $p(\xi_t;\rho_{0:t})$. To do so, we propagate the uncertainty of the affine prior policy (\ref{eq:SPPDP_Cont_prior}) through the nonlinear function of the dynamical system (\ref{eq:general_dynamics_transition}) as follows
\begin{subequations}
    \begin{align}
        x_{t+1} &= f_t(x_t,\hat{u}_t+\delta u_t) \\
        \delta u_t &\sim \rho_t(\delta u_t|\delta x_t) \approx \mathcal{N}(\delta u_t;k_t+\matrixstyle{K}_t \delta x_t,\Sigma_t)
    \end{align}
\end{subequations}

Equivalently, we construct the following recursion
\begin{subequations} \label{eq:SPPDP_forward_closedloop}
    \begin{align}
        p(\xi_t;\rho_{0:t}) &= p(u_t|x_t;\rho_t) p(x_t;\rho_{0:t-1})\\
        p(x_{t+1}|\rho_{0:t}) &= \int p(x_{t+1}|\xi_t) p(\xi_t;\rho_{0:t})d\xi_t
    \end{align}
\end{subequations}

We propose to approximate this density using a multivariate Gaussian, that is $p(\xi_t|\rho_{0:t}) \approx \mathcal{N}(\xi_t;\mu_{\xi,t},\Sigma_{\xi\xi,t})$, by applying the unscented transform. 
The recursion is initialized with some arbitrary initial state density, $p(x_0)$. We also assume this to be a Gaussian $\mathcal{N}(x_0|\mu_{x,0},\Sigma_{xx,0})$. Then the procedure in (\ref{eq:SPPDP_forward_closedloop}) can be specialised to
\begin{subequations}
\label{eq:SPPDP_Forward}
    \begin{align}
        \mu_{u,t} &= \hat{u}_t + k_t+\matrixstyle{K}_t (\mu_{x,t}-\hat{x}_t)\\
        \Sigma_{uu,t} &= \matrixstyle{K}_t \Sigma_{xx,t} \matrixstyle{K}_t^\top +\Sigma_t \\
        \Sigma_{ux,t} &= \matrixstyle{K}_t \Sigma_{xx,t} \\
        f_{t,n} &\triangleq f_t\left(\mu_{\xi,t}+\sqrt{\Sigma}_{\xi\xi,t}\epsilon^{\xi}_n\right) \\
        \mu_{x,t+1} &= \sum\nolimits_{n=1}^{N_\xi} w^{\xi}_n f_{t,n} \\
        \Sigma_{xx,t+1} &= \sum\nolimits_{n=1}^{N_\xi} w^{\xi}_n (f_{t,n} - \mu_{x,t+1})(f_{t,n} - \mu_{x,t+1})^\top
    \end{align}
\end{subequations}
where 
\begin{subequations} \label{eq:XI_Statistics}
  \begin{align}
  	 \mu_{\xi,t} &\triangleq \begin{bmatrix}
  		\mu_{x,t}  \\
  		\mu_{u,t}  \\
  	\end{bmatrix} \\
  	  \Sigma_{\xi\xi,t} &\triangleq \begin{bmatrix}              				\Sigma_{xx,t} & \Sigma_{ux,t}^{\top}  \\
             		* & \Sigma_{uu,t} \\
                			    \end{bmatrix}
  \end{align}
\end{subequations}

The weights $\{w^{\xi}_n\}_{n=1}^{N_\xi}$ and associated unit sigma points $\{\epsilon^{\xi}_n\}_{n=1}^{N_\xi}$ correspond with some predetermined signature. 
\subsubsection{Backward Pass}
In the backward pass, the key idea is to approximate the probabilistic state-action value function using a second-order Fourier-Hermite series expansion about the nominal closed-loop trajectory density. Together with the affine Gaussian prior policy approximation, this implies that also the optimal cost-to-go is quadratic. Finally, this allows us to evaluate the backward recursion from Theorem \ref{th:IterativeProbabilisticDP}.

The following Lemma evaluates both the quadratic optimal cost-to-go and the control affine optimal policies. 
\begin{lemma} \label{lemma:SPPDP_Update}
    Assume quadratic approximations for the probabilistic state-action value and optimal cost-to-go functions
   \begin{subequations}
   	 \begin{align}\label{eq:Quadratization_Q}
            Q_t^*(\xi_t) &\approx \hat{Q}_t^*(\xi_t)=\frac{1}{2} \begin{bmatrix}
				1 \\ \delta \xi_t
			\end{bmatrix}^\top \begin{bmatrix}
				2\hat{Q}^*_{0,t} & \hat{Q}_{\xi,t}^{*,\top}  \\
				* & \hat{Q}_{\xi\xi,t}^*  \\
			\end{bmatrix}\begin{bmatrix}
				1 \\ \delta \xi_t
		\end{bmatrix}\\
	\label{eq:Quadratization_V}
            V_t^*(x_t) &\approx \hat{V}_t^*(x_t)=\frac{1}{2} \begin{bmatrix}
				1 \\ \delta x_t
			\end{bmatrix}^\top \begin{bmatrix}
				2\hat{V}^*_{0,t} & \hat{V}_{x,t}^{*,\top}  \\
				* & \hat{V}_{xx,t}^*  \\
			\end{bmatrix}\begin{bmatrix}
				1 \\ \delta x_t
		\end{bmatrix}
    \end{align} 
   \end{subequations}
    where $\delta \xi_t \triangleq \xi_t - \mu_{\xi,t}$. 
    
    Now, also assume that we adopt the affine Gaussian prior policy in (\ref{eq:SPPDP_Cont_prior}). Then, $\pi_t^*$ is also affine Gaussian and its gains can be evaluated as
    \begin{subequations}\label{eq:SPPDP_ContUpdate}
    	\begin{align}
    		\matrixstyle{K}_t^* &= \Sigma_t^*(\Sigma_t^{-1}\matrixstyle{K}_t - \lambda\hat{Q}^*_{ux,t}) \\
    		{k}_t^* &= \Sigma_t^*(\Sigma_t^{-1}{k}_t - \lambda\hat{Q}^*_{u,t}) \\
    		\Sigma_t^* &= ( \Sigma_t^{-1} + \lambda\hat{Q}^*_{uu,t})^{-1} 
    	\end{align}
    \end{subequations}
    Furthermore, $\hat{V}_t^*$, is quadratic and its coefficients can be evaluated as
    \begin{subequations}\label{eq:SPPDP_ValueUpdate}
    	\begin{align} 
                \hat{V}_{0,t}^* &= \hat{Q}_{0,t}^* + \tfrac{1}{\lambda}  \left(\tfrac{1}{2}k_t^{\top} \Sigma_t^{-1} {k}_t -\tfrac{1}{2}k_t^{*\top} \Sigma_t^{*,-1} {k}_t^* {\color{black}+ \log  \tfrac{\det(\Sigma_t)}{\det(\Sigma_t^*)} }\right) \\
                \hat{V}_{x,t}^* &= \hat{Q}_{x,t}^* + \tfrac{1}{\lambda} \left(\matrixstyle{K}_t^{\top} \Sigma_t^{-1} {k}_t - \matrixstyle{K}_t^{*\top} \Sigma_t^{*,-1} {k}_t^*\right) \\
                \hat{V}_{xx,t}^* &= \hat{Q}_{xx,t}^* + \tfrac{1}{\lambda} \left(\matrixstyle{K}_t^{\top} \Sigma_t^{-1} \matrixstyle{K}_t - \matrixstyle{K}_t^{*\top} \Sigma_t^{*,-1} \matrixstyle{K}_t^*\right)
    	\end{align}
    \end{subequations}
\end{lemma}
\begin{proof}
    By substituting (\ref{eq:SPPDP_Cont_prior}), (\ref{eq:Quadratization_Q}), and (\ref{eq:Quadratization_V}) into (\ref{eq:PDP_ValueFunctions_V}) and (\ref{eq:PDP_Policy}), one can verify (\ref{eq:SPPDP_ValueUpdate}) and (\ref{eq:SPPDP_ContUpdate}), respectively.
\end{proof}

Now, as we have mentioned before, we propose to approximate $Q_t^*$ using a Fourier-Hermite series expansion. It is straightforward to choose the Gaussian measure as  $\mathcal{N}(\xi_t;\mu_{\xi,t},\Sigma_{\xi\xi,t})$ to indicate the region of interest. We further use the sigma-points method to evaluate the integrals. Here it is emphasized that evaluation of the state-action value function in (\ref{eq:PDP_ValueFunctions_Q}) entails evaluation of the state cost, $c_t$, and dynamic function, $f_t$. Interestingly, therefore we can recycle the sigma-points, $\{f_{t,n}\}$, that were used in the forward pass. 

It follows that the coefficients of (\ref{eq:Quadratization_Q}) are determined as
\begin{equation}\label{eq:FH_Q_2}
    \begin{aligned}
        \hat{Q}_{0,t}^* &= \tilde{Q}^*_{0,t} - \tfrac{1}{2}\trace\{\Tilde{Q}^*_{\xi\xi,t}\} \\
        \hat{Q}_{\xi,t}^* &= \sqrt{\Sigma}_{\xi\xi,t}^{-1}\tilde{Q}_{\xi,t}^*\\
        \hat{Q}_{\xi\xi,t}^* &= \sqrt{\Sigma}_{\xi\xi,t}^{-1}\tilde{Q}_{\xi\xi,t}^*\sqrt{\Sigma}_{\xi\xi,t}^{-1}
    \end{aligned}
\end{equation}
where
\begin{equation}\label{eq:FH_Q_3}
    \begin{aligned}
        \tilde{Q}^*_{0,t} &= \sum\nolimits_{n=1}^{N_\xi} w^{\xi}_n Q_t^*\left(\mu_{\xi,t}+\sqrt{\Sigma}_{\xi\xi,t} \epsilon^{\xi}_n\right)\\
        \tilde{Q}_{\xi,t}^* &= \sum\nolimits_{n=1}^{N_\xi} w^{\xi}_n Q_t^*\left(\mu_{\xi,t}+\sqrt{\Sigma}_{\xi\xi,t} \epsilon^{\xi}_n\right) \epsilon^{\xi}_n \\
        \tilde{Q}_{\xi\xi,t}^* &= \sum\nolimits_{n=1}^{N_\xi} w^{\xi}_n Q_t^*\left(\mu_{\xi,t}+\sqrt{\Sigma}_{\xi\xi,t} \epsilon^{\xi}_n\right) \left(\epsilon^{\xi}_n\epsilon^{\xi\top}_n -\matrixstyle{I}\right)
    \end{aligned}
\end{equation}

Once these coefficients are available, the optimal cost-to-go coefficients and optimal policy gains can be updated following Lemma \ref{lemma:SPPDP_Update} again and so forth. As such the recursion can be executed. That is except for $\hat{V}_T^*$. Therefore we can once again rely on the Fourier-Hermite series expansion. It follows
\begin{equation}
    \begin{aligned} \label{eq:Terminal_Value_Coef}
        \hat{V}_{0,T}^* &=  \tilde{V}_{0,T}^* - \tfrac{1}{2} \trace \{\tilde{V}_{xx,T}^*\} \\
        \hat{V}_{x,T}^* &= \sqrt{\Sigma}_{xx,T}^{-1} \tilde{V}_{x,T}^* \\
        \hat{V}_{xx,T}^* &= \sqrt{\Sigma}_{xx,T}^{-1}\tilde{V}_{xx,T}^{*}\sqrt{\Sigma}_{xx,T}^{-1}
    \end{aligned}
\end{equation}
where
\begin{equation}
    \begin{aligned} \label{eq:Terminal_Value_FH}
        \tilde{V}_{0,T}^* &= \sum\nolimits_{n=1}^{N_x} w^x_n c_T(\mu_{x,T}+\sqrt{\Sigma}_{xx,T} \epsilon^x_n)\\
        \tilde{V}_{x,T}^* &= \sum\nolimits_{n=1}^{N_x} w^x_n  c_T(\mu_{x,T}+\sqrt{\Sigma}_{xx,T}\epsilon^x_n) \epsilon^x_n\\
        \tilde{V}_{xx,T}^* &= \sum\nolimits_{n=1}^{N_x} w^x_n  c_T(\mu_{x,T}+\sqrt{\Sigma}_{xx,T} \epsilon^x_n) (\epsilon^x_n {\epsilon^{x\top}_n} -\matrixstyle{I})
    \end{aligned}
\end{equation}

In conclusion, we must still evaluate the optimal risk parameter, $\gamma^*$. Such is established by the following Lemma.

\begin{lemma}
   The optimal risk parameter (\ref{eq:Opt_Parameter}) can be approximated by  
    \begin{equation} \label{eq:opt_gamma_SPPDP}
            \begin{aligned}
                \gamma^* & \approx  \frac{T+1}{\sum_{t=0}^T \sum_{n=1}^{N_\xi} w^{\xi}_n c_t\left(\sqrt{\Sigma}_{\xi \xi,t} \epsilon^{\xi}_n + \mu_{\xi,t}\right)}
            \end{aligned}
    \end{equation}
    where the signature corresponds with the signature produced by the forward pass in the next iteration.
\end{lemma}
\begin{proof}
    The denominator in $\gamma^*$ in (\ref{eq:Opt_Parameter}) can be rephrased
    \begin{equation*}
        \sum\nolimits_{t=0}^{T} \mathbb{E}_{p(\xi_{t}|\mathcal{O}_{0:T}; \rho_{0:T-1},\lambda)}\left[ c_t(\xi_t) \right]
    \end{equation*}
    Note that the expectation is evaluated with respect to $p(\xi_{t}|\mathcal{O}_{0:T}; \rho_{0:T-1},\lambda)$ rather than $p(\xi_{t}; \rho_{0:T-1},\lambda)$ so that we cannot recycle the signature generated by the current forward pass. However, by construction, we have that 
    \begin{equation}
    	p(\xi_{t}|\mathcal{O}_{0:T};\rho_{0:T-1},\lambda) \equiv p(\xi_{t};\pi^*_{0:T-1})
    \end{equation}
    
    This then proves (\ref{eq:opt_gamma_SPPDP}). We can calculate $\gamma^*$ after the forward pass. 
    
\end{proof}

The entire procedure is summarized in Algorithm \ref{Alg:SPPDP}. We refer to this algorithm as the Sigma-Point Probabilistic Dynamic Programming or SP-PDP algorithm.

\begin{algorithm}[t]
	\caption{SP-PDP}
	\label{Alg:SPPDP}
	\KwIn{ $T$, $f_{0:T-1}$, $c_{0:T-1}$, $c_T$, $\mu_{x,0}$, $\Sigma_{xx,0}$, $\rho_{0:T-1}$, $\lambda$, $\{\epsilon^{\xi}_n, w^{\xi}_n\}_{n=1}^{N_\xi}$ and $\{\epsilon^{x}_n, w^{x}_n\}_{n=1}^{N_x}$ 
	} 
    \KwOut{$u^{*}_{0:T-1}$} 
    Initialize $\mu_{\xi,0:T}$ and $\Sigma_{\xi\xi,0:T}$ using (\ref{eq:SPPDP_Forward}). \\
	\Repeat{convergence}{
		Evaluate $\hat{V}_{0,T}^*$, $\hat{V}_{x,T}^*$ and $\hat{V}_{xx,T}^*$ using (\ref{eq:Terminal_Value_Coef}).\\
		\For(\tcp*[f]{Backward Pass}){$t = T-1:0$}
		{
			Evaluate $\hat{Q}_{0,t}^*$, $\hat{Q}_{\xi,t}^*$ and $\hat{Q}_{\xi \xi,t}^*$ using (\ref{eq:FH_Q_2})-(\ref{eq:FH_Q_3}).\\
			Evaluate $ k_t^*, \matrixstyle{K}_t^*$ and $\Sigma^*_t$ using (\ref{eq:SPPDP_ContUpdate}).\\ 
			Evaluate $\hat{V}_{t}^*$, $\hat{V}_{x,t}^*$ and $\hat{V}_{xx,t}^*$ using (\ref{eq:SPPDP_ValueUpdate}).\\
		}
		$\rho_{0:T-1} \leftarrow \pi^*_{0:T-1} $ \\
            Start from $\mu_{x,0}$ and $\Sigma_{xx,0}$.\\
		\For(\tcp*[f]{Forward Pass}){$t = 0:T-1$}
		{
			Evaluate $\mu_{\xi,t}$, $\Sigma_{\xi \xi,t}$, $\mu_{x,T}$ and $\Sigma_{xx,T}$ using (\ref{eq:SPPDP_Forward}). \\
			Store $\{f_{t,n}\}$ to recycle in the backward pass.\\
		}
		Evaluate $\gamma^*$ using (\ref{eq:opt_gamma_SPPDP}) and $\lambda \leftarrow \gamma^*$. \\
	}
\end{algorithm}


\begin{remark}
	In the standard DDP algorithm, evaluation of the policy gains requires evaluation of the inverse of $\hat{Q}_{uu,t}$, see (\ref{eq:DDP_opt_gain}). Other DDP implementations therefore tend to suffer from situations where this matrix is not invertible which causes a singularity in the update \cite{Whittle1996,Hassan2023,jacobson1970differential}. Some researchers suggest using a regularization term to remedy this issue \cite{jacobson1970differential}. In the present approach such regularization is an intrinsic property of the gain evaluation, see (\ref{eq:SPPDP_ContUpdate}). 
\end{remark}

\begin{remark}
    The use of sigma-point methods and the unscented transform for the classic DDP has also been studied in \cite{Hassan2023, Manchester2016}. In these works, the Gaussian measure used to generate the sigma points is chosen arbitrarily. It was emphasized the performance of the algorithm strongly depends on this choice. The present approach avoids this issue by automatically determining the Gaussian measure. 
    
\end{remark}

\subsection{Sigma-Point Bayesian Smoothing Control}

Alternatively, we can now also pursue an algorithm from the policy evaluation procedure described in Theorem \ref{th:BayesianInferenceControl}. The main idea here deviates from the DDP algorithm in that we simply need to approximate the Bayesian smoothing equations, specifically the forward filtering equations and backward smoothing correction equations. Similar to the DDP and SP-PDP algorithms, the algorithm's main iteration will consist of a forward and a backward pass. However, other than merely evaluating the prior policy during the forward pass, the filtering equations will also encode information about the stage cost through the measurement model. Moreover, the backward smoothing correction does not encode information about the stage cost. This means that the algorithm simultaneously evaluates and explores the optimization landscape during the forward pass. The backward pass merely passes the information from future time steps back to the present time step.

\subsubsection{Forward Pass}
In the forward pass, we assume that the filtering densities can be reasonably approximated with a Gaussian distribution. Let us define
\begin{subequations}
	\begin{align}
        p(\xi_t|\mathcal{O}_{0:t-1}; \rho_{0:t},\lambda) &\approx \mathcal{N}\left(\xi_t;\mu_{\xi,t}^{f,-},\Sigma_{\xi\xi,t}^{f,-}\right) \label{eq:pred_dist}\\
        p(\xi_t|\mathcal{O}_{0:t}; \rho_{0:t},\lambda) &\approx \mathcal{N} \left(\xi_t;\mu_{\xi,t}^{f},\Sigma_{\xi\xi,t}^{f}\right)    \label{eq:filter_dist}
\end{align}
\end{subequations}
where $\mu_{\xi,t}^{f,-}$, $\mu_{\xi,t}^{f}$, $\Sigma_{\xi\xi,t}^{f,-}$, and $\Sigma_{\xi\xi,t}^{f}$ are defined in the same way as (\ref{eq:XI_Statistics}). In Bayesian filtering, they are referred to as the prediction and filtering densities, respectively. To evaluate the parameters, it is required that we evaluate the prediction and update steps (\ref{eq:BIC_UKF_pred}) and (\ref{eq:BIC_UKF_correct}).

The forward pass initialises with $\mathcal{N}(\xi_0;\mu_{\xi,0}^{f},\Sigma_{\xi\xi,0}^{f})$ where
\begin{subequations} \label{eq:SPBIC_Forward_Pred_Initial}
	\begin{align}  
				\mu_{\xi,0}^{f} &= \begin{bmatrix}
			\mu_{x,0}  \\
			k_{0}+\matrixstyle{K}_{0} \mu_{x,0}  \\
		\end{bmatrix} \\
		\Sigma_{\xi\xi,0}^{f} &= \begin{bmatrix}
			\Sigma_{xx,0} &  \Sigma_{xx,0} \matrixstyle{K}_{0}^{\top} \\
* & \matrixstyle{K}_{0} \Sigma_{xx,0} \matrixstyle{K}_{0}^\top +\Sigma_{0} \\
		\end{bmatrix} 
	\end{align}
\end{subequations}

Next, we can approximate the prediction step (\ref{eq:BIC_UKF_pred}). We employ the unscented transform to propagate the uncertainty through the nonlinear map, $f_t$. One verifies that
\begin{subequations} \label{eq:SPBIC_Forward_Pred}
    \begin{align} 
        f_{t-1,n} &= f_{t-1}\left(\mu_{\xi,t-1}^{f}+\sqrt{\Sigma}_{\xi\xi,t-1}^f\epsilon^{\xi}_n\right)\\
        \mu_{x,t}^{f,-} &= \sum\nolimits_{n=1}^{N_\xi} w^{\xi}_n f_{t-1,n} \\
        \Sigma_{xx,t}^{f,-} &= \sum\nolimits_{n=1}^{N_\xi} w^{\xi}_n \left(f_{t-1,n} - \mu_{x,t}^{f,-}\right)\left(f_{t-1,n} - \mu_{x,t}^{f,-}\right)^\top \\
        \mu_{u,t}^{f,-} &= k_{t}+\matrixstyle{K}_{t} \mu_{x,t}^{f,-}\\
        \Sigma_{uu,t}^{f,-} &= \matrixstyle{K}_{t} \Sigma_{xx,t}^{f,-} \matrixstyle{K}_{t}^\top +\Sigma_{t} \\
        \Sigma_{ux,t}^{f,-} &= \matrixstyle{K}_{t} \Sigma_{xx,t}^{f,-}
    \end{align}
\end{subequations}


{
	
	Now, to approximate the update step in (\ref{eq:BIC_UKF_correct}). Inspired by the SP-PDP algorithm we propose to approximate the stage cost, $c_t$, using a quadratic surrogate
	\begin{equation}
		\label{eq:noraml_dist_cost}
		 c_t(\xi_t) \approx \hat{c}_t(\xi_t)=\frac{1}{2} \begin{bmatrix}
			1 \\ \delta \xi_t
		\end{bmatrix}^\top \begin{bmatrix}
			2\hat{c}_{0,t} & \hat{C}_{\xi,t}^{\top}  \\
			* & \hat{C}_{\xi\xi,t}  \\
		\end{bmatrix}\begin{bmatrix}
			1 \\ \delta \xi_t
		\end{bmatrix}\\
	\end{equation}
	
	The coefficients of $\hat{c}_t$ can be determined using a second-order Fourier-Hermite series expansion. Here it is straightforward to use as a Gaussian measure the prediction density, $\mathcal{N}(\xi_t;\mu_{\xi,t}^{f,-},\Sigma_{\xi\xi,t}^{f,-})$. The coefficients can then be evaluated similarly to (\ref{eq:FH_Q_2}) and (\ref{eq:FH_Q_3}).

\begin{remark}
	Note that as opposed to the SP-PDP algorithm, the dynamic function, $f_t$, and, cost rate, $c_t$, are therefore not evaluated for the same signature. 
\end{remark}

}


Deriving from the quadratic stage cost approximation in (\ref{eq:noraml_dist_cost}), the following lemma allows us to evaluate the filtering density (\ref{eq:filter_dist}).

\begin{lemma}
    The parameters of the filtering density (\ref{eq:filter_dist}) are given by
    \begin{subequations}\label{eq:SPBIC_Forward_Corr}
        \begin{align}
            \mu_{\xi,t}^{f} &=  \left(I-\Gamma_t\right) \left(\mu_{\xi,t}^{f,-}- \lambda \Sigma_{\xi \xi , t}^{f-} \hat{C}_{\xi,t}\right) \\
            \Sigma_{\xi \xi , t}^{f} &= \left(I-\Gamma_t\right) \Sigma_{\xi \xi , t}^{f-}
        \end{align}
    \end{subequations}
    where $\Gamma_t$ is referred to as a filtering gain and is defined as
    \begin{equation}\label{eq:Filtering_Gain}
        \Gamma_t \triangleq \Sigma_{\xi \xi , t}^{f-} \left(\Sigma_{\xi \xi , t}^{f-} + \tfrac{1}{\lambda} \hat{C}_{\xi \xi,t}^{-1}\right)^{-1} 
    \end{equation}
\end{lemma}
\begin{proof}
   We approximate the measurement model defined in (\ref{eq:exp_cost}) as a Gaussian using the quadratic model in (\ref{eq:noraml_dist_cost}).
    \begin{equation*}
    	p(\mathcal{O}_{t}|\xi_{t};\lambda) \approx  \mathcal{N}\left(\xi_t;-\hat{C}_{\xi \xi,t}^{-{\top}}\hat{C}_{\xi,t},\tfrac{1}{\lambda}\hat{C}_{\xi \xi,t}^{-1}\right)
    \end{equation*}
    
    Substituting the Gaussian approximation above and the Gaussian prediction density, $\mathcal{N}(\xi_t;\mu_{\xi,t}^{f,-},\Sigma_{\xi\xi,t}^{f,-})$, in the nominator of (\ref{eq:BIC_UKF_correct}) 
    yields a product of two Gaussian densities. This is proportional to another Gaussian density $\mathcal{N}(\xi_t;m_{\xi,t}, P_{\xi \xi,t})$. 
    
    Its mean and covariance are as follows
    \begin{equation*}
        \begin{aligned}
            m_{\xi,t} &= P_{\xi \xi,t} \left(\Sigma_{\xi\xi,t}^{f,-,-1}\mu_{\xi,t}^{f,-} -\lambda\hat{C}_{\xi,t}\right)\\
            P_{\xi \xi,t} &= \left(\Sigma_{\xi\xi,t}^{f,-,-1}+\lambda\hat{C}_{\xi\xi,t}\right)^{-1}
        \end{aligned}
    \end{equation*}
    
	It directly follows that
    \begin{equation*}
        \begin{aligned}
            p(\xi_t|\mathcal{O}_{0:t}; \rho_{0:t},\lambda) &\approx \mathcal{N} \left(\xi_t;\mu_{\xi,t}^{f},\Sigma_{\xi\xi,t}^{f}\right)= \mathcal{N}\left(\xi_t;m_{\xi,t},P_{\xi \xi,t}\right)  
        \end{aligned}
    \end{equation*}
    
    Finally using the matrix inversion lemma for the covariance matrix, $P_{\xi \xi,t}$, one obtains
    \begin{equation*}
    \Sigma_{\xi\xi,t}^f = \Sigma_{\xi \xi , t}^{f-}-\Sigma_{\xi \xi , t}^{f-} \left(\Sigma_{\xi \xi , t}^{f-} + \lambda \hat{C}_{\xi \xi, t}^{-1}\right)^{-1} \Sigma_{\xi \xi , t}^{f-}
    \end{equation*}
    
    For computational stability, we introduce the filtering gain in (\ref{eq:Filtering_Gain}) and rewrite the obtained mean and covariance matrix as mentioned in (\ref{eq:SPBIC_Forward_Corr}). 
\end{proof}


\subsubsection{Backward Pass}
In the backward pass, again we assume that the smoothing density is reasonably approximated with a Gaussian. Let us define
\begin{equation} \label{eq:smooth_dist}
        p(\xi_{t}|\mathcal{O}_{0:T}; \rho_{0:T-1},\lambda) \approx \mathcal{N} (\xi_t;\mu_{\xi,t}^{s},\Sigma_{\xi\xi,t}^{s})
\end{equation}
where $\mu_{\xi,t}^{s}$ and $\Sigma_{\xi \xi,t}^{s}$ are defined in the same way as presented in (\ref{eq:XI_Statistics}).

To evaluate the smoothing density, we follow the same procedure as in \cite{sarkka2023}. Note that (\ref{eq:BIC_URTSS1}) can be obtained by marginalizing the following joint probability 
\begin{equation} \label{eq:SPBIC_smooth_joint}
    \begin{multlined}
        p(\xi_t, x_{t+1}|\mathcal{O}_{0:T}; \rho_{0:T-1},\lambda) = \\p(\xi_t|x_{t+1}, \mathcal{O}_{0:t}; \rho_{0:t},\lambda) p(x_{t+1}|\mathcal{O}_{0:T}; \rho_{0:T-1},\lambda),
    \end{multlined}
\end{equation}

The density, $p(\xi_t|x_{t+1}, \mathcal{O}_{0:t}; \rho_{0:t},\lambda)$, is approximated by another Gaussian
\begin{equation}
	\label{eq:smoothing_prediction_dens}
        p(\xi_t, x_{t+1}|\mathcal{O}_{0:t}; \rho_{0:t}) \approx 
            \mathcal{N} \left(\begin{bmatrix}
            	\xi_t \\ x_{t+1}
            \end{bmatrix}; \mu^{s,-}_t, \Sigma^{s,-}_t \right)
\end{equation}
where 
\begin{subequations} \label{eq:smooth_pred_dist_param}
    \begin{align}
    	\mu_t^{s,-} &\triangleq \begin{bmatrix}
    		\mu_{\xi,t}^{f}  \\
    		\mu_{x,t+1}^{f,-}  \\
    	\end{bmatrix} \\
    	    	\Sigma_t^{s,-} &\triangleq \begin{bmatrix}
    		\Sigma_{\xi \xi,t}^{f} & \Sigma_{\xi x, t+1}^{s,-}   \\
    		* & \Sigma_{xx,t+1}^{f,-}
    	\end{bmatrix}
    \end{align}
\end{subequations}

The parameter, $\Sigma_{xx,t+1}^{f,-}$, can be obtained by recycling the sigma points from the forward pass.
    \begin{equation} \label{eq:SPBIC_Backward_Pred}
        \Sigma_{\xi x, t+1}^{s,-} = \sqrt{\Sigma}_{\xi\xi,t}^f \sum\nolimits_{n=1}^{N_\xi} w^{\xi}_n \epsilon^{\xi}_n \left(f_{t,n} - \mu_{x,t+1}^{f,-}\right)^\top
    \end{equation}

Then the parameters of the smoothing density can be evaluated by substituting (\ref{eq:smooth_dist}) and (\ref{eq:smoothing_prediction_dens}) into (\ref{eq:SPBIC_smooth_joint}) and marginalizing out $x_{t+1}$. One ultimately verifies that
\begin{equation}
        \begin{aligned} \label{eq:SPBIC_Backward_Corr}
            \mu_{\xi,t}^s &= \mu_{\xi,t}^f + \matrixstyle{S}_t(\mu_{x,t+1}^{s}-\mu_{x,t+1}^{f,-})\\
            \Sigma_{\xi \xi,t}^s &= \Sigma_{\xi \xi,t}^{f} + \matrixstyle{S}_t(\Sigma_{xx,t+1}^{s}-\Sigma_{xx,t+1}^{f,-})\matrixstyle{S}_t^\top
        \end{aligned}
\end{equation}
where the $\matrixstyle{S}_t$ is referred to the smoothing gain and is defined as
\begin{equation} \label{eq:SPBIC_Backward_Gain}
        \matrixstyle{S}_t \triangleq \Sigma_{\xi x, t+1}^{s,-} \Sigma_{xx,t+1}^{f,-,-1}
\end{equation}
Further, note that the backward recursion in (\ref{eq:SPBIC_Backward_Corr}) is initialised from $\mu_{x,T}^{s} \triangleq \mu_{x,T}^{f}$ and $\Sigma_{xx,T}^{s} \triangleq \Sigma_{xx,T}^{f}$. 

Based on the results above we are now capable of evaluating the optimal policies and the optimal risk parameter. Such is established by the following Lemma. 
This result was also given in \cite{watson2021b} for stochastic systems dynamics and excluding the generalised cost formulation that is used here.
\begin{lemma}
    The parameters of the optimal controller (\ref{eq:SPPDP_Cont_posterior}) can be evaluated by 
    \begin{equation} \label{eq:control_params}
            \begin{aligned}
                \matrixstyle{K}_t^* &= \Sigma_{ux,t}^{s} \Sigma_{xx,t}^{s,-1} \\
                k^* &= \mu_{u,t}^{s}-\matrixstyle{K}_t^* \mu_{x,t}^{s} \\
                \Sigma^*_t &= \Sigma_{uu,t}^{s} - \matrixstyle{K}_t^* \Sigma_{xx,t}^{s}  {\matrixstyle{K}_t^*}^{\top}
            \end{aligned}
    \end{equation}
    and the optimal risk parameter can be evaluated as
    {\begin{equation} \label{eq:opt_gamma_SPBIC}
            \begin{aligned}
                \gamma^* &=  \frac{T +1}{\sum_{t=0}^T \tfrac{1}{2}\trace(\Sigma_{\xi \xi,t}^{s}\hat{C}_{\xi \xi,t}) + \tfrac{1}{2} \mu_{\xi,t}^{s^{\top}}\hat{C}_{\xi \xi,t}\mu_{\xi,t}^s + \hat{C}_{\xi,t}^{\top} \mu_{\xi,t}^s + \hat{c}_{t}}
            \end{aligned}
    \end{equation}}
\end{lemma}
\begin{proof}
    Following Theorem \ref{th:EM_solution}, we can calculate the conditional probability $p(u_t|x_t,\mathcal{O}_{0:T};\rho_{0:T-1},\lambda)$ using the smoothing density (\ref{eq:smooth_dist}). One verifies that
    \begin{equation*}
            \begin{multlined}
                p(u_t|x_t,\mathcal{O}_{0:T};\rho_{0:T-1},\lambda) = \\ 
                \mathcal{N} \Big(u_t; \mu_{u,t}^{s} + \Sigma_{ux,t}^{s} \Sigma_{xx,t}^{s,-1} (x_t-\mu_{x,t}^{s}), \\ \Sigma_{uu,t}^{s} - \Sigma_{ux,t}^{s} \Sigma_{xx,t}^{s,-1} \Sigma_{xu,t}^{s} \Big)
            \end{multlined} 
    \end{equation*}
    
    Second, again rephrasing the denominator of $\gamma^*$ in (\ref{eq:Opt_Parameter}) as
    \begin{equation*}
    	\sum\nolimits_{t=0}^{T} \mathbb{E}_{p(\xi_{t}|\mathcal{O}_{0:T}; \rho_{0:T-1},\lambda)}\left[ c_t(\xi_t) \right]
    \end{equation*}
	and now substituting (\ref{eq:noraml_dist_cost}) and (\ref{eq:smooth_dist}) one retrieves (\ref{eq:opt_gamma_SPBIC}).
\end{proof}

The entire procedure is summarized in Algorithm \ref{Alg:SPBSC}. We refer to this algorithm as the Sigma-Point Bayesian Smoothing Control or SP-BSC algorithm.
\begin{algorithm}[t]
\caption{SP-BSC}
\label{Alg:SPBSC}
	\KwIn{ $T$, $f_{0:T-1}$, $c_{0:T-1}$, $c_T$, $\mu_{x,0}$, $\Sigma_{xx,0}$, $\rho_{0:T-1}$, $\lambda$,$\{\epsilon^{\xi}_n , w^{\xi}_n\}_{n=1}^{N_\xi}$ 
	}
	\KwOut{ $k^*_{0:T-1}$ and $\matrixstyle{K}^*_{0:T-1}$} 
\Repeat{convergence}{
    Evaluate $\mu_{\xi,0}^{f}$ and $\Sigma_{\xi\xi,0}^{f}$ using (\ref{eq:SPBIC_Forward_Pred_Initial}).\\
    \For(\tcp*[f]{Forward Pass}){$t = 1:T$}
        {
            Evaluate $\mu_{\xi,t}^{f,-}$, $\Sigma_{\xi \xi,t}^{f,-}$, $\mu_{x,T}^{f,-}$ and $\Sigma_{xx,T}^{f,-}$ using (\ref{eq:SPBIC_Forward_Pred}). \\
            Evaluate $\mu_{\xi,t}^f$, $\Sigma_{\xi \xi,t}^f$, $\mu_{x,T}^f$ and $\Sigma_{xx,T}^f$ using (\ref{eq:SPBIC_Forward_Corr}).\\            Store $f_{t,n}$, $\mu_{\xi,t}^{f,-}$, $\Sigma_{\xi \xi,t}^{f,-}$, $\mu_{x,T}^{f,-}$ and $\Sigma_{xx,T}^{f,-}$ to recycle in the backward pass.\\
        }
    Set $\mu_{x,T}^{s} = \mu_{x,T}^{f}$ and $\Sigma_{xx,T}^{s} = \Sigma_{xx,T}^{f}$.\\
    \For(\tcp*[f]{Backward Pass}){$t = T-1:0$}
        {
            Evaluate $\Sigma_{\xi x, t+1}^{s,-}$ using (\ref{eq:SPBIC_Backward_Pred}).\\
            Evaluate $\mu_{\xi,t}^{s}$ and $\Sigma_{\xi \xi,t}^{s}$ using (\ref{eq:SPBIC_Backward_Corr}).\\ 
            Evaluate $k_t^*,\matrixstyle{K}_t^*$ and $  \Sigma^*_t$ using (\ref{eq:control_params}).\\ 
        }
    Evaluate $\gamma^*$ using (\ref{eq:opt_gamma_SPBIC}).\\
    $\rho_{0:T-1}, \lambda \leftarrow  \pi^*_{0:T-1}, \gamma^*$ \\
}
\end{algorithm}

\begin{remark}
    The use of sigma-point methods, such as Gauss-Hermite and cubature quadrature, for implementing the message-passing technique has also been explored in \cite{watson2021b} within stochastic dynamics.
\end{remark}
    \subsection{Discussion and connection to prior work}\label{sec:prev_work}
	
    The first algorithm is closest related to the DDP algorithm and makes a clear distinction between updating the policy (backward pass) and evaluating the policy (forward pass). The second algorithm differs from the first algorithm in the sense that evaluation and updating of the policy happen simultaneously (forward pass), followed by a second update step that now also backtraces information from the entire trajectory (backward pass), similar to the first algorithm. This means that during the evaluation step, the new policy (and associated optimal trajectory) are already being updated based on the cost accumulated up to time, $t$. This is as opposed to the first algorithm. This property can be considered as a feature in the sense that this might improve the speed of convergence for some problems. On the other hand, the trajectory distribution that is produced during the forward step does no longer reflect the real trajectory distribution when the last fully optimized policy would be used on the physical system. 
		
	Based on the discussed theoretical and practical features, we hypothesize that both algorithms will strike a better balance between exploration and exploitation than their predecessors. The uncertainty on the policies will automatically scale the locality of the sigma-point methods that are used to evaluate the pseudo-gradient. The larger the uncertainty, the larger the area searched in the optimization space. It is proposed that this feature helps the algorithms navigate the optimization landscape more efficiently. Secondly, the framework automatically incorporates regularization with respect to the prior policy. We expect that this will aid with numerical stability.

{
In conclusion, we highlight the relationship with prior work. First we remark that the algorithmic structure of the second approach (SP-BSC) was proposed prior to our work by \cite{watson2021advancing,watson2021b} using and comparing Gauss-Hermite and cubature quadrature uncertainty quantification methods. Compared with this work, we remark that we explicitly limit the application of our algorithm to a deterministic setting, given the ambiguity associated with optimizing the risk-sensitivity parameter in a stochastic setting. Second, we use a different uncertainty-quantification strategy seeking analogies with the recent work in Fourier-Hermite DDP \cite{Hassan2023}. Further, in \cite{watson2020stochastic}, the authors compare the algorithm with Guided Policy Search (GPS) which includes a variant of iLQR tailored to the entropy regularized SOC problem \cite{levine2013guided,levine2013variational}. In that sense, our first approach (SP-PDP) is structurally closely related to the GPS framework. However, as mentioned in the introduction, iteration of the entropy regularized SOC, and thus iteration of the iLQR tailored to the entropy regularized SOC problem, will not converge to the RSOC policy but rather to the SOC policy. This means that we cannot directly compare results in a stochastic setting. The difference can be traced back to a difference in the definition of the value function in the SOC and RSOC setting. This difference could be implemented in the GPS framework but would deviate from the algorithm described in \cite{levine2013guided,levine2013variational} and the algorithm would then converge to the RSOC policy. Furthermore, these works do not discuss optimization of the risk parameter nor is it straightforward to do so in the theoretical setting provided there. The theoretical analysis given here clearly shows how optimization of the risk parameter can be applied in general, regardless of how the policy is evaluated in the end, i.e. by dynamic programming or smoothing. Ultimately this allows for a fairer comparison. 

In summary, our contribution is to demonstrate how the second approach is theoretically equivalent to the first approach and how both solve a RSOC problem when applied in a stochastic setting and without optimizing the risk parameter. In a deterministic setting, the risk parameter can be used as a hyper-parameter. Second, when approximations are used to cast these approaches into computational algorithms, the difference manifests. The analysis from section \ref{sec:expectation-maximization} helps to interpret and understand their connection. Our formal theoretical treatment also allows to set the risk parameter with either algorithm which would otherwise not be straightforward.} 


\section{Simulation Experiments} \label{sec:results}

To demonstrate the performance of the algorithms proposed in section \ref{sec:algorithm}, we conducted experiments on three nonlinear Hamiltonian systems. We compare the presented methods with the SP-DP presented in \cite{Hassan2023}. This approach implements the DDP method using sigma-point methods to approximate the state-action value function instead of a Taylor series expansion and outperforms the standard DDP method. All of our experiments were conducted using MATLAB R2021b on a 1.80 GHz Intel Core i7-1265U CPU. 

In the following, we generate swing-up manoeuvres for a single pendulum and a cart pole and generate optimal motion planning for a 6-DoF industrial robot. In all our experiments, we used a quadratic cost model of the form
\begin{subequations}\label{eq:quad_cost}
	\begin{align}
		c_t(\xi_t) &= (x_t-x_g)^{\top} W (x_t-x_g) + u_t^{\top} R u_t\\
		c_T(x_T) &= (x_T-x_g)^{\top} W_T (x_T-x_g)
	\end{align}
\end{subequations}
where $x_g$ represents the goal state we aim to reach and $W$, $W_T \in \mathbb{R}^{n_x \times n_x}$ and $R \in \mathbb{R}^{n_u \times n_u}$ represent weight matrices.

To generate the unit sigma points and their corresponding weights, we employed the third and fifth-order unscented transforms (UT3 and UT5) using the spherical cubature rule and Gauss–Hermite (GH) quadrature rule of order $p=3$. For details, we refer to \cite{sarkka2023, Hassan2023}. To discretise the continuous dynamics, we used a fourth-order Runge–Kutta integration method and a zero-order hold for the control action, $u$. 

\subsection{Pendulum Swing-up}
This experiment aims to swing up the single pendulum from the downward position $(\theta = 0)$ to the horizontal position $(\theta = \pi/2)$ by applying an input torque $u$. It follows that $n_x=2$ and $n_u=1$. The state of this system is defined as $x=(\theta,\dot{\theta})$, where $\theta$ and $\dot{\theta}$ denote the angle and angular velocity of the pendulum, respectively. The model parameters are adopted from \cite{Hassan2023, Manchester2016}
We set $T=50$ with sampling time $0.1 \unit{\second}$. The initial and goal states are $x_0 = (0,0)$ and $x_g = (\pi/2,0)$, respectively. We use the weight matrices $W=10^{-1} I$, $R=10^{-1}I$, and $W_T = 10^2 I$ in which $I$ indicates an identity matrix with proper dimension. 
For SP-PDP and SP-BSC, We assumed that $\lambda = 1$ for initialization, the initial prior policy is $\mathcal{N}(u_t; 0, 10^{-1}I)$, and $\Sigma_{xx,0}=10^{-3}I$. As mentioned we compare the presented methods with the SP-DP algorithm. With this method the covariance matrix, $\Sigma_{\xi \xi, t}$, is set to an arbitrary but fixed value. We set $\Sigma_{\xi \xi, t} = 10^{-6} I$ which was the best choice documented in \cite{Hassan2023}.

\begin{figure}[t]
	\captionsetup{font=footnotesize}
	\centering
	\includegraphics[width=1\columnwidth]{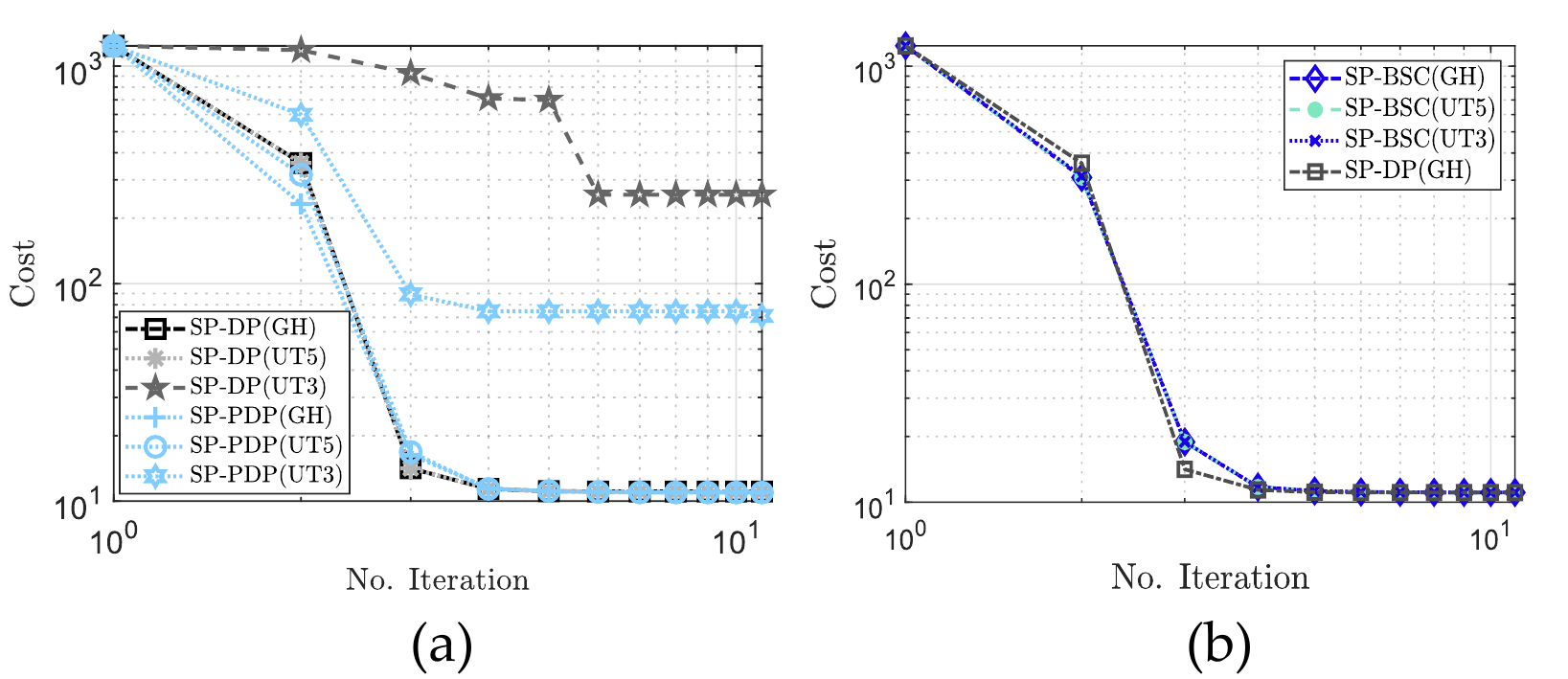}
	\caption{\small{Total cost per iteration for the pendulum swing-up experiment. In (a), the results of SP-PDP are compared to the SP-DP using the UT3, UT5, and GH cubature rules. In (b), the results of SP-BSC using the UT3, UT5, and GH cubature rules are compared with the best results obtained with SP-DP.}}\label{fig:pend_cost}
\end{figure}

Fig. \ref{fig:pend_cost} shows the total cost of the optimized trajectory in each iteration for different algorithms. 
Fig. \ref{fig:pend_cost} (a) demonstrates the results for SP-PDP in comparison to the previous algorithm SP-DP using the different sigma-point generation rules UT3, UT5, and GH. In Fig  \ref{fig:pend_cost} (b), the results for SP-BSC using UT3, UT5, and GH rules are compared to the best result obtained with SP-DP. 
As can be seen in Fig. \ref{fig:pend_cost}, the SP-DP and SP-PDP methods combined with the UT3 rule do not find the optimal solution however the SP-BSC method does. 
We conclude that SP-PDP and SP-DP require higher-order sigma-point generation rules to be efficient. 
Further remark that all methods except SP-DP and SP-PDP with UT3 rule approximately demonstrate share a similar convergence history although SP-PDP and SP-BSC achieve a higher reduction in cost than SP-DP in the first iteration.

\begin{figure}[t]
	\captionsetup{font=footnotesize}
	\centering
	\includegraphics[width=1\columnwidth]{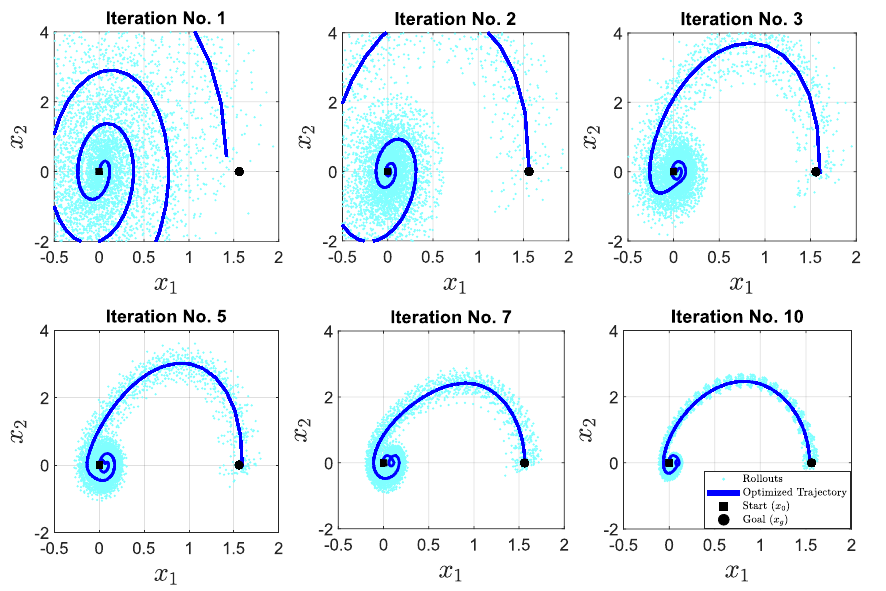}
	\caption{\small{The optimized trajectory and the sampled points of 60 trajectory rollouts in the state-space at each iteration of the SP-PDP using the GH cubature rule.}}\label{fig:pend_evolution}
\end{figure}

Fig. \ref{fig:pend_evolution} presents the obtained trajectory in the state space at each iteration of the SP-PDP using the GH cubature rule. For illustrational purposes, we showed the sampled points of 60 trajectory rollouts generated by the optimized probabilistic policy at each iteration. These trajectories are generated by simulating the deterministic system dynamics in a closed loop whilst sampling from the probabilistic policies. As shown in Fig. \ref{fig:pend_evolution}, the sampled points are spread out in a large part of the state-space in the initial iterations and gradually the achieved trajectory converges to the optimum and the sampled points concentrate around the optimal trajectory. 
The dispersion of the points in this figure is associated with the covariance matrix $\Sigma_{\xi \xi, t}$ that indicates the interesting area around the nominal points for expanding the FH series. 
The covariance matrix $\Sigma_{\xi \xi, t}$ is chosen arbitrarily in SP-DP. In contrast, in the discussed paradigm, this covariance matrix stems from the policy's uncertainty, which is optimized at each iteration. Consequently, during early iterations, the covariance matrix is large to facilitate exploration of the state space as shown in Fig. \ref{fig:pend_evolution}. As the probabilistic policies converge to the deterministic optimal policy, its confidence increases and the covariance matrix shrinks correspondingly as shown in Fig. \ref{fig:pend_evolution}. We conclude that the algorithms discussed here are capable of maintaining the right balance between exploration and exploitation. 

Table \ref{tab:pend} reports the number of used sigma points and the average run times (in milliseconds) to evaluate the backward and forward passes per iteration for all methods. The number of sigma points used at step time $t$ and at terminal time step $T$ to evaluate terminal value function is denoted by $N_\xi$ and $N_x$, respectively. 
In general, the computational speed for all algorithms mainly depends on the number of sigma points used in the integration rule.
Based on Table \ref{tab:pend}, for a specific sigma-point generation rule, SP-BSC is faster than other algorithms and the run time for SP-BSC using the UT3 rule is the fastest among all the methods. The computational time of SP-PDP is slightly larger than SP-DP because, in the forward pass of SP-PDP, we need more computations than SP-DP to propagate the uncertainty. 

In conclusion, reasoning from Fig. \ref{fig:pend_cost} and Table. \ref{tab:pend}, we can conclude that SP-BSC showcases the best overall performance.

\begin{table}[b]
	\captionsetup{font=footnotesize}
	\caption{Average run times as a function of the number of sigma-points for different algorithms in the pendulum swing-up experiment.}
	\label{tab:pend}
	\centering
	\begin{tabular}{l|c|c}\hline \hline
		\bfseries Algorithm & \bfseries  Average run time $[\unit{\milli\second}]$ & \bfseries  No. of Sigma-points\\ \hline 
		SP-DP (UT3) &   8.6 & $N_\xi=7$, $N_x=5$\\ 
		SP-DP (UT5) &   12.0  & $N_\xi=19$, $N_x=9$\\ 
		SP-DP (GH)  &   13.6  & $N_\xi=27$, $N_x=9$\\ 
		SP-PDP (UT3)&   10.9 & $N_\xi=7$, $N_x=5$\\ 
		SP-PDP (UT5)&   12.5 & $N_\xi=19$, $N_x=9$\\ 
		SP-PDP (GH) &   15.4 & $N_\xi=27$, $N_x=9$\\ 
		SP-BSC (UT3)&   6.2  & $N_\xi=7$ \\
		SP-BSC (UT5)&   9.5  & $N_\xi=19$ \\
		SP-BSC (GH) &   11.2  & $N_\xi=27$ \\
		\hline \hline
	\end{tabular}
\end{table}

\subsection{Cart-pole Swing-up}

In this experiment, the goal is to swing up the cart-pole system by applying a horizontal force $u$. For this experiment, we have $n_x=4$ and $n_u=1$.  The state vector of the system is defined as $x=(p, \theta, \dot{p}, \dot{\theta})$, where $p$ and $\dot{p}$ represent the position and the velocity of the cart, respectively, and $\theta$ and $\dot{\theta}$ denote the angle and angular velocity of the pole, respectively. Again we adopt the model parameters from \cite{Hassan2023, Manchester2016}. We set $T=30$ with sampling time $0.1s$. The initial and goal states are $x_0 = (-1,\pi/2,0,0)$ and $x_g = (0,\pi,0,0)$, respectively.
The cost weight matrices $W$, $W_T$, $R$, and $\Sigma_{xx,0}$ are set the same as with the pendulum swing-up experiment. 
Further, we assumed that $\lambda = 1$ with initialization and the initial prior policy is $\mathcal{N}(u_t; 0, 5I)$ for SP-PDP and SP-BSC.
We chose $\Sigma_{\xi \xi, t} = 10^{-6} I$ for SP-DP. Again this was the best choice mentioned in \cite{Hassan2023}.

\begin{figure}[t]
	\captionsetup{font=footnotesize}
	\centering
	\includegraphics[width=1\columnwidth]{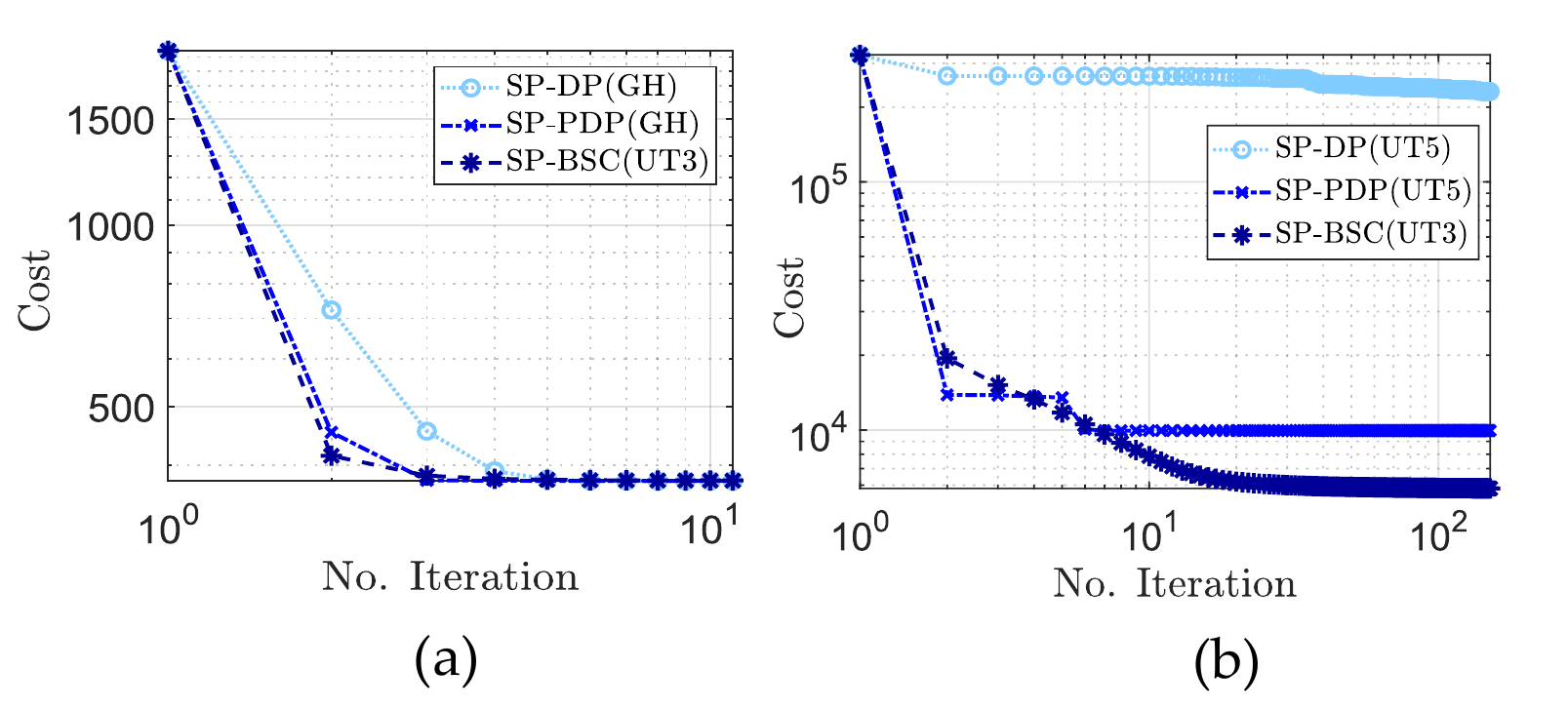}
	\caption{\small{Total cost per iteration for (a) the cart-pole swing-up experiment and (b) the 6-DoF robot motion planning experiment.}}\label{fig:cart_robot}
\end{figure}

Based on the results achieved with the pendulum swing-up experiment, we know that SP-DP and SP-PDP needed more sigma points to achieve better convergence. Therefore, we have implemented SP-DP and SP-PDP using the GH rule, and SP-BSC with the UT3 rule. 
In Fig. \ref{fig:cart_robot} (a) and Table \ref{tab:cart}, the results are documented for the cart-pole swing-up experiment. It can be seen that all algorithms eventually converged to a similar total cost but the cost reduction rates of the algorithms are more discriminative than in the pendulum swing-up experiment. SP-PDP and SP-BSC showcase a faster convergence rate than SP-DP. We conclude that the SP-PDP and SP-BSC can achieve better convergence rates than the SP-DP algorithm when the complexity and non-linearity of the problem increase.

Further note that the SP-BSC using the UT3 rule has a faster run time, as shown in Table \ref{tab:cart}. Based on Table \ref{tab:cart}, we notice that again the computation time of SP-PDP exceeds that of SP-DP due to the additional computations required in the forward pass of SP-PDP to propagate the uncertainty.

\begin{table}[b]
	\captionsetup{font=footnotesize}
	\caption{Average run times as a function of the number of sigma-points for different algorithms in the cart-pole swing-up experiment.}
	\label{tab:cart}
	\centering
	\begin{tabular}{l|c|c}\hline \hline
		\bfseries Algorithm & \bfseries  Average run time $[\unit{\milli\second}]$ & \bfseries  No. of Sigma-points\\ \hline 
		SP-DP (GH) &   63.0 & $N_\xi=243$, $N_x=81$\\ 
		SP-PDP (GH)&   94.1 & $N_\xi=243$, $N_x=81$\\ 
		SP-BSC (UT3)&  9.5  & $N_\xi=11$ \\
		\hline \hline
	\end{tabular}
\end{table}

\subsection{6-DoF Motion Planning}
Finally, we consider optimal motion planning of a 6-DoF industrial manipulator arm. It follows that $n_x=12$ and $n_u=6$.  The state vector of the system is $x=(\theta,\dot{\theta})$ in which $\theta \in \mathbb{R}^{6}$ and $\dot{\theta} \in \mathbb{R}^{6}$ denote the generalised coordinates and generalised velocities of the manipulator joints.
We aim to generate an optimal trajectory from the initial position $\theta_0 = 0$ to the goal position $\theta_g = (0, \pi/2, 0, -\pi/2, 0, \pi/2)$ by optimizing the torque input $u$. 
We set $T=200$ with time step $5 \unit{\milli\second}$. The cost weight matrices $W$, $W_T$, $R$, and $\Sigma_{xx,0}$ are again equivalent to the previous experiments. 
We assume $\lambda = 1$ for initialization and the initial prior policy is $\mathcal{N}(u_t; 0, 50I)$ for SP-PDP and SP-BSC.
We set $\Sigma_{\xi \xi, t} = 10^{-6} I$ for SP-DP again.

Using the GH rule is not feasible since the number of evaluation points for the GH rule is $p^{n_\xi}$, which increases exponentially as the number of dimensions increases. Since $n_\xi=18$ we have implemented the SP-DP and SP-PDP algorithms using the UT5 rule. The SP-BSC was implemented using the UT3 rule. 

The cost reduction curves of all algorithms for this experiment are shown in Fig. \ref{fig:cart_robot} (b). As mentioned in \cite{Hassan2023}, SP-DP can work with the high-order sigma-point generation rules but it appears not to work well for high-dimensional problems. This statement is confirmed by Fig. \ref{fig:cart_robot} (b). One notes that the SP-DP cannot converge with a similar number of iterations as the algorithms that are discussed in this work (SP-PDP and SP-BSC) while SP-BSC using the UT3 rule and SP-PDP using the UT5 rule can. After a few iterations, the cost reduction curve of SP-PDP is getting flat but SP-BSC can converge continuously and smoothly. 
Moreover, we listed the average run times (in seconds) to evaluate both backward pass and forward pass per iteration for all algorithms in Table \ref{tab:robot}. One observes that the SP-BSC is faster than the others too. The SP-PDP and SP-DP require almost the same computational time for execution. 

\begin{figure}[!h]
	\captionsetup{font=footnotesize}
	\centering
	\includegraphics[width=.85\columnwidth]{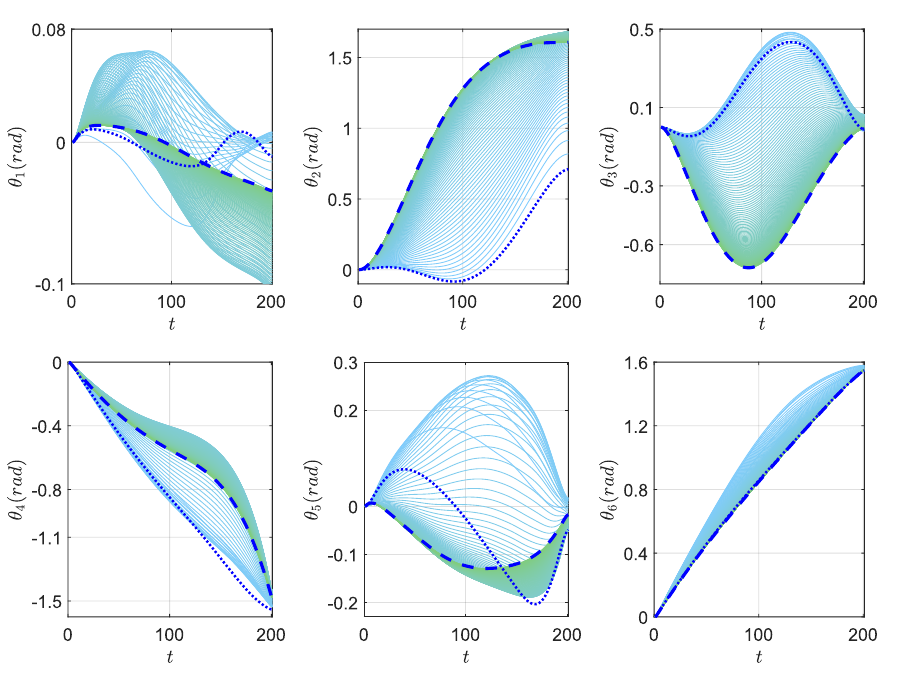}
	\caption{\small{Visualisation of the convergence of the robot joint angle trajectories with the SP-BSC method using the UT3 rule. The dotted, solid, and dashed lines indicate the first iteration, intermediate iterations, and the 150th iteration, respectively. }}\label{fig:robot_traj}
\end{figure}

In conclusion and for visual reference we illustrate the optimized robot joint angle trajectories in each iteration of SP-BSC using the UT3 rule in Fig. \ref{fig:robot_traj}. For this experiment, we aim to regularize the joint angles around the goal position $\theta_g$ and we set $W=10^{-1} I$, $R=10^{-2}I$, and $W_T = 10^{3} I$. 
The algorithm ran for 150 iterations. The first trajectory after applying the algorithm is shown by dotted lines, the trajectories of intermediate iterations are represented by solid lines, and the final dashed line depicts the optimal trajectory to regularize the robot around the goal angle $\theta_g$. 

\begin{table}[!h]
	\captionsetup{font=footnotesize}
	\caption{Average run times as a function of the number of sigma-points for different algorithms in the robot motion planning experiment.}
	\label{tab:robot}
	\centering
	\begin{tabular}{l|c|c}\hline \hline
		\bfseries Algorithm & \bfseries  Average run time $[\unit{\second}]$ & \bfseries  No. of Sigma-points\\ \hline 
		SP-DP (UT5) &   32.5034 & $N_\xi=649$, $N_x=289$\\ 
		SP-PDP (UT5)&   34.4003 & $N_\xi=649$, $N_x=289$\\ 
		SP-BSC (UT3)&   1.9133  & $N_\xi=37$ \\
		\hline \hline
	\end{tabular}
\end{table}

\section{Conclusion and Discussion} \label{sec:Conclusion}

In this paper, we have discuss two derivative-free algorithms tailored to TO, the Sigma-Point Probabilistic Dynamic Programming and the Sigma-Point Bayesian Smoothing Control algorithm. The algorithms are based on the contemporary probabilistic optimal control paradigm that re-established optimal control as an estimation problem. This allows us to lean on well-established inference tools to treat the problem which differ substantially from the conventional tools. 

The present approach has two properties that are direct consequences of the probabilistic reformulation and EM treatment.
\begin{itemize}
	\item 
	The proposed schemes are intrinsically iterative, whereas, in standard gradient-based optimization, an iterative calculation is typically introduced ad hoc. That is because the approaches discussed in this work relies on the EM method which introduces a principled iterative procedure by design. The EM approach further boasts an internal mechanism that regularizes the new solution to the previous one. This supports improved algorithmic stability.
	\item 
	The algorithms discussed in this work dynamically adjust the mechanism used to probe the optimization landscape inherently. 
	Specifically, the scale of exploration is directly related to the reliability of the current solution estimates. 
	This contrasts with standard gradient-based algorithms, which typically probe the landscape by simply evaluating gradient information at the current iterate. Due to the probabilistic reformulation, these method are associated to the risk-sensitive optimal control framework. 
	It is hypothesized that these approaches are risk-seeking in the initial iterations and converge to the original optimal control problem in the final iterations. Therewith they balance between exploration and exploitation.
\end{itemize}


These statements are supported by our numerical experiments where the algorithms discussed here outperform previous work, especially in complex and high-dimensional settings. 

Future efforts will focus on extending the proposed algorithms into online applications. Further, it would also be interesting to integrate constraints into the paradigm or at least account for practical constraints in the algorithms.

\section*{Acknowledgment}
This work was supported by the Research Foundation Flanders (FWO) under SBO grant no. S007723N. 
The authors wish to thank Joe Watson and Hany Abdulsamad for their insightful comments on earlier versions of this manuscript.

\bibliographystyle{IEEEtran}
\bibliography{References}

\vspace*{-35pt}

\begin{IEEEbiography}
[{\includegraphics[width=1in,height=1.25in,clip,keepaspectratio]{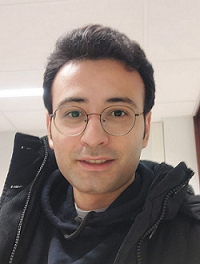}}]
{Mohammad Mahmoudi Filabadi} received his M.Sc. degree in Electrical Engineering with a major in Control Systems from the Sharif University of Technology, Tehran, Iran, in 2021.  

He is currently working toward a Ph.D. degree in Electromechanical Engineering at Ghent University, Ghent, Belgium. He is also an affiliate member of Flanders Make, Belgium. His current research interests include employing probabilistic inference techniques for controller design and trajectory optimization to achieve active exploration strategies, with a specific focus on applications within mechatronic and robotic systems. 
\end{IEEEbiography}

\vspace*{-35pt}

\begin{IEEEbiography}
[{\includegraphics[width=1in,height=1.25in,clip,keepaspectratio]{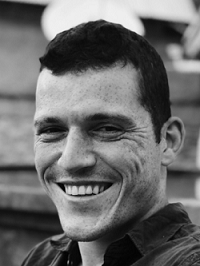}}]
{Tom Lefebvre} received the M.Sc. degree in control engineering and automation from Ghent University, Belgium, in 2015, and the Ph.D. degree in electromechanical engineering in 2019. 

Since 2019, he has been a post-doctoral researcher. He worked on stochastic optimal control, gradient-based and stochastic trajectory optimization, uncertainty quantification, and black-box optimization. His current research interests include probabilistic modelling, identification, and control. He is currently an affiliate member of Flanders Make.
\end{IEEEbiography}

\vspace*{-35pt}

\begin{IEEEbiography}
[{\includegraphics[width=1in,height=1.25in,clip,keepaspectratio]{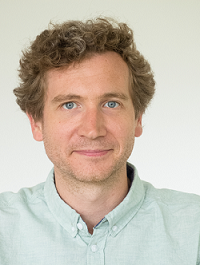}}]
{Guillaume Crevecoeur} (Member, IEEE) received his master's and Ph.D. degrees in Engineering Physics from Ghent University in 2004 and 2009, respectively. 

He received a Research Foundation Flanders postdoctoral fellowship in 2009, was appointed Associate Professor in 2014, and is a Full Professor since 2024 at Ghent University. He is a member of Flanders Make in which he leads the Ghent University activities on machines, intelligence, robotics and control. With his team, he conducts research on system identification and nonlinear control, mainly for mechatronics, industrial robotics and energy systems. His goal is to endow physical dynamic systems with improved functionalities and capabilities when interacting with uncertain environments, other systems and humans.
\end{IEEEbiography}

\end{document}